\documentclass[twocolumn,amsthm,secthm]{autart}
\usepackage{amssymb}
\usepackage{amsbsy}

\usepackage{graphicx}          % Include this line if your
\usepackage{subfigure}

\usepackage{color}

\newcommand{\rev}[1]{#1}
\newcommand{\revA}[1]{#1}

\begin{document}

\begin{frontmatter}
%\runtitle{Insert a suggested running title}  % Running title for regular
                                              % papers but only if the title
                                              % is over 5 words. Running title
                                              % is not shown in output.

% Title, preferably not more than 10 words.
\title{Intrinsic Reduced Attitude Formation with Ring Inter-Agent Graph\thanksref{footnoteinfo}}

\thanks[footnoteinfo]{This paper was not presented at any IFAC
meeting.
%Corresponding author W.~Song. Tel. +46 700967354, +86 15210909661.
%Fax +XXXIX-VI-mmmxxv.
}

% Add the e-mail address (ead) as shown
\author[kth]{Wenjun Song}\ead{wenjuns@kth.se},
\author[kth]{Johan Markdahl}\ead{markdahl@kth.se},
\author[kth,hit]{Silun Zhang}\ead{silunz@kth.se},
\author[kth]{Xiaoming Hu}\ead{hu@kth.se},
\author[amss]{Yiguang Hong}\ead{yghong@iss.ac.cn}

% Please supply full addresses here.
\address[kth]{Optimization and Systems Theory, Department of Mathematics, KTH Royal Institute of Technology, 100 44 Stockholm, Sweden}
\address[amss]{Key Laboratory of Systems and Control, Academy of Mathematics and Systems Science, Chinese Academy of Sciences, 100190 Beijing, China}
\address[hit]{Control and Simulation Center, Harbin Institute of Technology,150001 Harbin, China}

% Five to ten keywords, chosen from the IFAC keyword list or with the help of the Automatica keyword wizard

% Abstract of not more than 200 words.
\begin{abstract}
This paper investigates the reduced attitude formation control problem for a group of rigid-body agents \revA{using feedback based on relative attitude information.
Under both undirected and directed cycle graph topologies, it is shown that reversing the sign of a classic consensus protocol yields asymptotical convergence to formations whose
shape depends on the parity of the group size.}
Specifically, in the case of even parity the reduced attitudes converge asymptotically to a pair of antipodal points and distribute equidistantly on a great circle in the case of odd parity.
Moreover, when the inter-agent graph is an undirected ring, the desired formation is shown to be achieved from almost all initial states.
%Illustrative examples are provided to show the effectiveness of the proposed control law.
\end{abstract}

\end{frontmatter}

\section{Introduction}

Multi-agent coordination
\cite{beard2001,fax04tac,jadbabaie2003,olfatiSaber07p} has gained
increasing recognition and appreciation during the last decade.
Following many significant results on consensus, how to effectively
generate various multi-agent formations, patterns or subgroup
divisions has attracted much attention. Among the problems studied,
 attitude formation of multiple rigid-body agents is of key importance
 with wide potential
applications such as formation flying \cite{beard2001,scharf04} and
multi-camera surveillance \cite{tron14tac,wang13prl}.
Attitude synchronization or consensus, a special and simple formation pattern
of attitudes, has been widely studied
\cite{tron12cdc,tron13tac,sarlette09AC,thunberg14ac}.

The (full) attitude of a rigid-body agent can be represented by a
rotation matrix
that evolves on the Lie group $\mathcal{SO}(3)$.
However, many attitude control applications do not require all three degrees of freedom of the full attitude to be determined. \rev{In rigid-body pointing applications, for example a body-fixed camera, the solar panel of a satellite or an antenna need to point towards some desired direction, the rotation
about the pointing axis is irrelevant since this rotation does not
change the direction in which the agent points.} Moreover,
in under-actuated situations where the rigid-body is actuated by
only two independent control torques, for example due to the failure
of a third
actuator, the rotation about the unactuated axis is disregarded.
\revA{The reduced attitude provides the proper framework to deal with such a situation \cite{bullo95}.}
All these applications lead to a reduced attitude
control problem \cite{bullo95,nalin11cs,lee11cdc,mayhew10acc}, in which the
reduced attitude of two degrees of freedom is naturally identified
with a point on the 2-sphere
$\mathcal{S}^2$.

Consider the attitude formation problem for a system of $n$ rigid-body agents on the product manifold $\mathcal{SO}(3)^n$ or $(\mathcal{S}^2)^n$ under a continuous feedback control law based on relative attitude information.
It can be shown that consensus states are intrinsic equilibria of the closed-loop system
regardless of the topology of the inter-agent graph.
The work \cite{tron12cdc} and \cite{pereira2017family} achieve this attitude synchronization for full attitudes and reduced attitudes respectively.
However, due to the fact that  $\mathcal{SO}(3)^n$ and $(\mathcal{S}^2)^n$ are
compact manifolds without boundary \cite{bhat00scl},
continuous time-invariant feedback control also yields some other closed-loop
equilibria that vary with the graph topology.
These equilibria represent different attitude configurations of the
system, which may include a desired formation depending on the application.
A natural and interesting question thus arises: is it possible to
achieve a desired formation by imposing some suitable inter-agent graph
to the system and designing a feedback control with only relative attitude information that stabilizes the formation?

It is increasingly recognized that one of the important ideas in
multi-agent systems is to design simple distributed control algorithms with not only cooperative but also antagonistic interactions between neighboring agents.
Modulus consensus requires the moduli of all agent states to reach a common value but the agents may be separated into several {\itshape antagonistic subgroups}.
The simplest case concerning two antagonistic subgroups, {\it bipartite consensus}, models the inter-agent connection as a signed graph
\cite{proskurnikov14cdc,altafini,valcher2014scl}.
For a general modulus consensus case or even the extended {\it set surrounding} case, signed graphs are replaced by the graphs with complex weights expressed with complex adjacency matrices \cite{youchengLou15tac}.
These results demonstrate that antagonistic interactions are effective to generate new coordination or formation patterns, but the considered dynamics basically evolve in Euclidean spaces. Additionally, cooperative control of motion on the circle and sphere with both attractive and repulsive couplings are studied in \cite{olfati06} and \cite{li2014tac}, respectively, but the inter-agent graph is required to be undirected and complete.

By extending the coordination studies with antagonistic interactions discussed in Euclidean spaces to those on compact manifolds, this paper provides a partial but affirmative answer to the aforementioned question via investigating reduced attitude formation with both undirected and directed ring inter-agent graph.
In particular, we focus on the generation of attitude formation patterns using only relative attitude information of a group of rigid-body agents.
Compared to the full attitude formation problem, the reduced attitude formation is more intuitive and easier to visualize.

In this paper, a simple angular velocity control for reduced attitude formation is proposed on the basis of antagonistic interactions between neighboring agents.
Due to the geometry of the 2-sphere some interesting phenomena are
observed: the closed-loop system behaves differently under the proposed
distributed control when the parity of the total number of agents is different.
Specifically, the antipodal formation is achieved when the number is even,
and the cyclic formation is achieved when the number is odd.
It is shown that these two reduced attitude formations %or equilibria
are intrinsic in
the sense that they result from the geometry of the 2-sphere and the topology
of the inter-agent graph.
%Naturally, as the graph topology varies, other interesting intrinsic equilibria may arise and this will be a subject of our future study.
It is worthwhile to mention that, in addition to the simple control structure, another strength of the proposed method is that we do not need to have the desired formation given beforehand or the formation errors in the control, in contrast to most existing methods \cite{beard2001,song13ccc,taeyoungLee13acc}.

\rev{
Comparing the control protocol and the resulting formation in this paper with that of \cite{altafini,valcher2014scl} in Euclidean space, there are mainly two differences: (i) for neighboring agents with antagonistic interaction, only relative states of neighboring agents is utilized in reduced attitude control, while absolute position of neighboring agents is required in the control of \cite{altafini,valcher2014scl};
(ii) in the case of ring inter-agent graph, if the total number of agents is odd and all neighboring agents are antagonistic, cyclic reduced formation can be attained due to the geometry of the 2-sphere, while the positions of all agents in \cite{altafini,valcher2014scl} reach consensus at the origin because the graph is unbalanced.
}

The rest of the paper is organized as follows: in Section~\ref{sec:pre},
necessary preliminaries on the reduced attitude and the 2-sphere are
introduced. In Section~\ref{sec:formulation}, the reduced attitude formation problem
with the ring inter-agent graph is formulated and a distributed angular velocity control law is proposed.
The antipodal formation and cyclic formation of reduced attitudes are discussed in Section~\ref{sec:anfor} and Section~\ref{sec:odd:cyclic}, respectively.
Following that, illustrative examples are provided in Section~\ref{sec:simulation}, and the conclusions are given in Section~\ref{sec:conclusion}.

%%%%%%%%%%%%%%%%%%%%%%%%%%%%%%%%%%%%%%%%%%%%%%%%%%%%%%%%%%%%%%%%%%%%%%%%%%%%%%%%
\section{Notations and preliminaries}\label{sec:pre}

%\subsection{Reduced Attitude}
%\subsection{The 2-sphere}

This paper considers the reduced attitude control problem for a
network of $n$ $(n\geq 2)$ rigid-body agents. In this section, we give
some preliminaries on the reduced attitude and the 2-sphere.

Let the index set $\mathcal{V}=\{1,2,\ldots,n\}$ represent the agents in the network.
Denote $R_i\in\mathcal{SO}(3)$ as the attitude of agent $i\in\mathcal{V}$ relative to the inertial frame $\mathcal{F}$, where $\mathcal{SO}(3)=\{R\in\mathbb{R}^{3\times 3} : R^TR=I, \det(R)=1\}$ is the rotation group of $\mathbb{R}^3$.
The kinematics of $R_i$ is governed by \cite{murray}
\begin{equation}\label{eq:kine:full}
\dot{R}_i = \widehat{\omega}_i R_i,
%\quad i\in\mathcal{V},
\end{equation}
where $\omega_i\in\mathbb{R}^3$ is the angular velocity of agent $i$ in the inertial frame $\mathcal{F}$, and the {\it hat} operator $(\cdot)^{\wedge}$ is defined by the equality that $\widehat{x}y=x\times y$ for any $x,y\in\mathbb{R}^3$.

Suppose that $b_i\in\mathcal{S}^2$ is a constant pointing direction in the
body-fixed frame of agent $i$, where $\mathcal{S}^2 = \{
x\in\mathbb{R}^3 : \|x\|=1 \}$ is the 2-sphere and $\|\cdot\|$ is the
Euclidean norm.
Let $\Gamma_i\in \mathcal{S}^2$ denote the same pointing direction
resolved in the inertial frame $\mathcal{F}$, then,
$\Gamma_i=R_ib_i$. $\Gamma_i$ is referred to as the {\it reduced attitude} of
agent $i$ since the rotation of agent $i$ about $b_i$ is ignored.
In the paper, we use a parametrization of $\Gamma_i$ given as follows
\begin{equation}\label{eq:gamma:para}
\Gamma_i=
%[\cos(\psi_i)\cos(\theta_i), \sin(\psi_i)\cos(\theta_i), \sin(\theta_i)]^T.
\left[
\begin{array}{c}
\cos(\psi_i)\cos(\varphi_i) \nonumber\\
\sin(\psi_i)\cos(\varphi_i)\nonumber\\
\sin(\varphi_i)\nonumber
\end{array}
\right].
\end{equation}
\rev{where $\varphi_i \in [-\pi/2,\pi/2]$ and $\psi_i\in[-\pi,\pi)$.
In fact, when $b_i=[1,0,0]^T$, the two angles $-\varphi_i$ and $\psi_i$ are the respective pitch and yaw angles of the rotation $R_i$.}
By the kinematics (\ref{eq:kine:full}) of the full attitude, the kinematics of $\Gamma_i$ is governed by \cite{lee11cdc}
\begin{equation}\label{eq:kine:reduced}
\dot{\Gamma}_i = \widehat{\omega}_i \Gamma_i.
%,\quad i\in\mathcal{V}.
\end{equation}

The tangent space of $\mathcal{S}^2$ at a point $\Gamma\in\mathcal{S}^2$ is given by $T_{\Gamma}\mathcal{S}^2=\{ x\in\mathbb{R}^3 : x^T \Gamma=0\}$.
Rotating $\Gamma\in\mathcal{S}^2$ about a unit axis $u\in T_{\Gamma}\mathcal{S}^2$ through an arbitrary angle $\beta$ transforms it to another point $\exp(\beta\widehat{u})\Gamma \in \mathcal{S}^2$, where $\exp(\cdot)$ is the matrix exponential.
For any two reduced attitudes $\Gamma_i,\Gamma_j\in\mathcal{S}^2$, define $\theta_{ij}\in[0,\pi]$ and $k_{ij}\in\mathcal{S}^2$ as
\[
\theta_{ij} = \arccos(\Gamma_i^T \Gamma_j)
,\quad
k_{ij}=\widehat{\Gamma}_i \Gamma_j / \sin(\theta_{ij}).
\]
It holds that %$k_{12}\in T_{\Gamma_1}\mathcal{S}^2 $ %\cap T_{\Gamma_2}\mathcal{S}^2 and
$\Gamma_j = \exp(\theta_{ij}\widehat{k}_{ij})\Gamma_i$.
Notice that the above equation for the unit axis $k_{ij}$ is valid only when $\theta_{ij}\in(0,\pi)$.
When $\theta_{ij}=0$ or $\pi$, we stipulate that $k_{ij}$ is chosen as any unit vector orthogonal to $\Gamma_i$. %$\sin(\theta_{ij}) k_{ij}=\widehat{\Gamma}_i\Gamma_j=0_3$.

The geodesic distance between any two points $\Gamma_i, \Gamma_j\in\mathcal{S}^2$, denoted as $d_{\mathcal{S}^2} (\Gamma_i,\Gamma_j)$, is the length of the shorter arc on the great circle of $\mathcal{S}^2$ joining the two points. Therefore,
\[
d_{\mathcal{S}^2} (\Gamma_i,\Gamma_j) = \theta_{ij}.
\]
The following lemma gives the relationship among geodesic distances of any three points on the 2-sphere,
which can be verified using spherical cosine formula.
More details about the geometry of the 2-sphere can be found in \cite{ferreira,todhunter}.

\begin{lem}\label{lem:sphere}
For any three points $\Gamma_i,\Gamma_j,\Gamma_k\in\mathcal{S}^2$,
\[
\cos(\theta_{ij})
=
\cos(\theta_{ik})\cos(\theta_{jk})
+ \sin(\theta_{ik})\sin(\theta_{jk})k_{ik}^Tk_{jk},
\]
\[
\theta_{ij} + \theta_{ik} + \theta_{jk} \leq 2\pi.
\]
Furthermore, $\Gamma_i,\Gamma_j,\Gamma_k$ lie on a great circle of $\mathcal{S}^2$ if and only if $\theta_{ij}=|\theta_{ik} - \theta_{jk}|$, $\theta_{ij}=\theta_{ik} + \theta_{jk}$ or $\theta_{ij} + \theta_{ik} + \theta_{jk}=2\pi$.
\end{lem}

Denote the state space of the system (\ref{eq:kine:reduced}) as the product manifold $(\mathcal{S}^2)^n$, which is the $n$-fold Cartesian product of $\mathcal{S}^2$ with itself.
We use $\mathbf{\Gamma}=\{\Gamma_i\}_{i\in\mathcal{V}} \in (\mathcal{S}^2)^n$ to denote the state of the system, and use the metric in $(\mathcal{S}^2)^n$ as
\[
d_{(\mathcal{S}^2)^n} (\mathbf{\Gamma},\bar{\mathbf{\Gamma}})
= \max_{i\in\mathcal{V}} d_{\mathcal{S}^2}(\Gamma_i,\bar{\Gamma}_i),
\quad
\forall \mathbf{\Gamma}, \bar{\mathbf{\Gamma}}\in(\mathcal{S}^2)^n.
\]
For a set $\mathcal{M}\subset (\mathcal{S}^2)^n$, define the distance from a point $\mathbf{\Gamma}\in (\mathcal{S}^2)^n$ to $\mathcal{M}$ as
\[
d_{\mathcal{M}}(\mathbf{\Gamma})
%\| \mathbf{\Gamma}\|_{\mathcal{M}}
= \inf_{\bar{\mathbf{\Gamma}}\in\mathcal{M} }
d_{(\mathcal{S}^2)^n} (\mathbf{\Gamma},\bar{\mathbf{\Gamma}}).
\]

\section{Problem formulation and control design}\label{sec:formulation}

In this section we provide a detailed description of the system and state the problem formulation. The qualitative behaviour of the system depends on the parity of the number of agents whereby the even and odd cases are treated separately.
The second half of this section concerns the control design.

\subsection{Intrinsic reduced attitude formation}

We model information exchange between the agents by a graph $\mathcal{G} = (\mathcal{V},\mathcal{E})$, where the vertex set $\mathcal{V}$ is given in Section~\ref{sec:pre} and the edge set $\mathcal{E}\subset \mathcal{V} \times \mathcal{V}$.
%Graph
$\mathcal{G}$ is undirected if the node pair of each edge is unordered,
%for any $(i,j)\in \mathcal{E}$, $(j,i)\in\mathcal{E}$;
otherwise it is directed.
A node $j$ is said to be a neighbor of $i$ if $(j,i)\in\mathcal{E}$, and $\mathcal{N}_i=\{ j: (j,i)\in\mathcal{E}\}$ is denoted as the set of neighbors of node $i$.
Throughout the paper, modulo $n$ operation is used to identify agents, for example, agent $0$ is the same as agent $n$ and agent $n+1$ is the same as agent $1$.
%we regard the index $n+1$ the same as $1$, and $0$ the same as $n$.

\begin{assum}\label{assump:graph}
The inter-agent graph $\mathcal{G}$ is either an undirected ring (or cycle), i.e., $\mathcal{N}_i=\{i-1,i+1\}$ for any $i\in\mathcal{V}$; or a directed ring (or cycle), i.e., $\mathcal{N}_i=\{i+1\}$ for any $i\in\mathcal{V}$.
\end{assum}

Under Assumption~\ref{assump:graph}, we define $W: (\mathcal{S}^2)^{n} \rightarrow\mathbb{R}$ as
\[
W(\mathbf{\Gamma}) = \min_{i\in\mathcal{V}} d_{\mathcal{S}^2}(\Gamma_i,\Gamma_{i+1}),
%\mathbf{\Gamma} \in (\mathcal{S}^2)^n.
\]
which is the minimal geodesic distance between the reduced attitudes of neighboring agents.
Before giving the main description of the reduced attitude formation problem in the paper, we give the following two lemmas about the maximal value of $W$, which result from the geometry of the 2-sphere.
The proofs can be found in Appendix~\ref{appen:sec:proof}.

\begin{lem}\label{lem:basic:even}
Suppose $n$ is even, then
\[
\max_{\mathbf{\Gamma}\in(\mathcal{S}^2)^n} W(\mathbf{\Gamma})
= \pi.
\]
Furthermore, $W(\mathbf{\Gamma})=\pi$ if and only if $\mathbf{\Gamma}\in\mathcal{M}_e$, where
\[
\mathcal{M}_e = \{ \mathbf{\Gamma}\in(\mathcal{S}^2)^n:
    \Gamma_i = (-1)^{i-1}v, \forall i\in\mathcal{V}, v\in\mathcal{S}^2\}.
\]
\end{lem}

\begin{lem}\label{lem:basic:odd}
Suppose $n$ is odd, then
\[
\max_{\mathbf{\Gamma}\in(\mathcal{S}^2)^n} W(\mathbf{\Gamma})
= \pi-\pi/n.
\]
It holds that $W(\mathbf{\Gamma}) = \pi-\pi/n$ if and only if $\mathbf{\Gamma}\in\mathcal{M}_o$, where
\begin{eqnarray}
\mathcal{M}_o = \{ \mathbf{\Gamma}\in(\mathcal{S}^2)^n &:&
    \Gamma_i = \exp\left((i-1)(\pi-\pi/n)\widehat{u}\right) v,
\nonumber \\
&&\forall i\in\mathcal{V},  u,v\in\mathcal{S}^2, u^Tv=0\}. \nonumber
\end{eqnarray}
\end{lem}

Notice that the manifold $\mathcal{M}_e$ (or $\mathcal{M}_o$) is empty when $n$ is odd (or even).

When $n$ is even and $\mathbf{\Gamma}\in\mathcal{M}_e$, the agents are divided into two groups indexed by $\{1,3,\ldots,n-1\}$ and $\{2,4,\ldots,n\}$, respectively.
The reduced attitudes are identical in the same group, and opposite in different groups, which is a pair of antipodal points on $\mathcal{S}^2$.
Under Assumption~\ref{assump:graph}, the geodesic distance between the reduced attitudes of any two neighboring agents is $\pi$.
We refer to this configuration for the system as an {\itshape antipodal formation} of reduced attitudes, a similar terminology as {\itshape bipartite consensus} \cite{altafini} of positions and {\itshape bipolar synchrony} \cite{olfati06} on the unit circle $\mathcal{S}^1$.

When $n$ is odd and $\mathbf{\Gamma}\in \mathcal{M}_o$, the reduced attitudes of the system lie equidistantly on a great circle with agents indexed by $\{1,3,5,\ldots,n\}$ on one half circle and agents indexed by $\{n-1,\ldots,6,4,2\}$ on the other.
The minimal geodesic distance between the reduced attitudes of any agent and every other agent is $2\pi/n$.
Under Assumption~\ref{assump:graph}, the geodesic distance between the reduced attitudes of any two neighboring agents is $\pi-\pi/n$.
We refer to this configuration for the system as a {\itshape cyclic formation} of reduced attitudes.

We are now ready to describe the reduced attitude formation problem that is addressed in this paper.

%{\it Definition \refstepcounter{mylem}\themylem.}
{\it Reduced attitude formation problem}: construct an angular
velocity controller for each agent $i$ with $
\{\widehat{\Gamma}_i\Gamma_j:j\in\mathcal{N}_i\}$ such that $\mathcal{M}_e$ whenever $n$ is
even and $\mathcal{M}_o$ whenever $n$ is odd are asymptotically
stable.
%as time tends to infinity, $\|\mathbf{\Gamma}\|_{\mathcal{M}_e}\rightarrow 0$ when $n$ is even,
%and $\|\mathbf{\Gamma}\|_{\mathcal{M}_o}\rightarrow 0$ when $n$ is odd.

\begin{rem}
Since the coordinates of $\widehat{\Gamma}_i\Gamma_j$ resolved in the body-frame of agent $i$ is $R_i^T \widehat{\Gamma}_i\Gamma_j = b_i \times (R_i^TR_j b_j)$, where $R_i^TR_j b_j$ is the pointing direction of agent $j$ seen from agent $i$, it is a relative reduced attitude information.
\end{rem}

\subsection{Control design and closed-loop system}\label{sec:controller}

In the paper, we propose the following distributed control rule:
\begin{equation}\label{eq:controller}
\omega_i = - \sum_{j\in\mathcal{N}_i} \widehat{\Gamma}_i \Gamma_j,
\quad i=1,2,\ldots,n.
\end{equation}
By the reduced attitude kinematics (\ref{eq:kine:reduced}), the closed-loop system can be written as
\begin{equation}\label{eq:closedloop}
\dot{\Gamma}_i = \widehat{\Gamma}_i \sum_{j\in\mathcal{N}_i} \widehat{\Gamma}_i \Gamma_j,
\quad i=1,2,\ldots,n.
\end{equation}
The next lemma states that the form of (\ref{eq:closedloop}) \revA{is invariant to any rotational coordinate transformation.}

\begin{lem}\label{lem:invariant}
For any $\alpha\in[0,\pi]$ and $u\in\mathcal{S}^2$, the closed-loop system (\ref{eq:closedloop}) is invariant under the coordinate transformation
\[
\bar{\Gamma}_i = \exp(\alpha\widehat{u})\Gamma_i, \quad i=1,2,\ldots,n,
\]
i.e., the closed-loop system in terms of $\bar{\mathbf{\Gamma}}=\{\bar{\Gamma}_i\}_{i\in\mathcal{V}}$ is
\[
\dot{\bar{\Gamma}}_i = \widehat{\bar{\Gamma}}_i \sum_{j\in\mathcal{N}_i} \widehat{\bar{\Gamma}}_i \bar{\Gamma}_j,
\quad i=1,2,\ldots,n.
\]
\end{lem}

Using the parametrization (\ref{eq:gamma:para}), the closed-loop system (\ref{eq:closedloop}) can also be written in terms of $\{(\psi_i,\varphi_i)\}_{i\in\mathcal{V}}$ as
\begin{eqnarray}
\cos(\varphi_i) \dot{\psi}_i
&=& - \sum_{j\in\mathcal{N}_i} \sin(\psi_j - \psi_i) \cos(\varphi_j),
\label{eq:closedloop:psi}\\
\dot{\varphi}_i
&=& \sum_{j\in\mathcal{N}_i} \big[ \sin(\varphi_i)\cos(\varphi_j)\cos(\psi_i-\psi_j)
\nonumber \\
&&\qquad\;\; -\cos(\varphi_i)\sin(\varphi_j)\big],
\label{eq:closedloop:varphi}
\end{eqnarray}
where $i=1,2,\ldots,n$.

\begin{rem}
Let the local potential function for agent $i\in\mathcal{V}$ be
\[
\varphi_i(\Gamma_i) = \frac{1}{2} \sum_{j\in\mathcal{N}_i} \|\Gamma_i-\Gamma_j\|^2 .
\]
The gradient of $\varphi_i$  on $\mathcal{S}^2$ at $\Gamma_i$ is given as
\[
\nabla_{\Gamma_i} \varphi_i = \widehat{\Gamma}_i \sum_{j\in\mathcal{N}_i} \widehat{\Gamma}_i\Gamma_j.
\]
Therefore, the angular velocity control (\ref{eq:controller}) makes the reduced attitude of agent $i$ move toward the positive gradient direction of $\varphi_i$ at $\Gamma_i$,
which means every agent tries to keep its reduced attitude away from the reduced attitudes of its neighboring agents.
%The next two sections will show that the reduced attitudes of the system will approach to $\mathcal{M}_e$ when $n$ is even or $\mathcal{M}_o$ when $n$ is odd, provided that Assumption \ref{assump:graph} holds.
\end{rem}

\begin{rem}\label{rem:consensus}
If we reverse the sign in (\ref{eq:controller}) such that
\begin{equation}\label{eq:controller:consensus}
\omega_i = \sum_{j\in\mathcal{N}_i} \widehat{\Gamma}_i \Gamma_j,
\quad i=1,2,\ldots,n,
\end{equation}
then with certain assumption on the connectivity of the inter-agent graph, $\Gamma_1, \ldots,\Gamma_n$ reach consensus asymptotically provided that their initial values lie on the surface of an open hemisphere, which is half $\mathcal{S}^2$. Proof of this conclusion can refer to \cite{tron13tac} for undirected graph and \cite{thunberg14ac} for directed and time-varying graph.
By imposing a reshaping function on the control algorithm, \cite{tron12cdc} shows that the new control rule makes the consensus manifold almost global attractive.
Notice that with the control law (\ref{eq:controller:consensus}), any reduced attitudes configuration that lies on a great circle (including $\mathcal{M}_e$ and $\mathcal{M}_o$) except the consensus case is unstable since this kind of configuration is on the boundary of the region of attraction of the consensus manifold.
\end{rem}

The consensus manifold of reduced attitudes is an intrinsic equilibrium set of the closed-loop system (\ref{eq:closedloop}), and by Remark~\ref{rem:consensus}, it is unstable.
The next two sections will show that the antipodal formation manifold $\mathcal{M}_e$ with even number of agents and the cyclic formation manifold $\mathcal{M}_o$ with odd number of agents are asymptotically stable under Assumption~\ref{assump:graph}.
Since the desired formation is not used in the control (\ref{eq:controller}), we refer to this kind of formation as {\it intrinsic reduced attitude formation}.

\begin{rem}
Unlike the reduced attitude coordination,
concerning the full attitude coordination problem with regular inter-agent graph, simulation demonstrates that a configuration of the full attitudes that lie equidistantly on a closed geodesic of $\mathcal{SO}(3)$ is asymptotically stable by a basic synchronization control law without any need to reverse the sign \cite{tron12cdc}.
\end{rem}

%We show that $\mathcal{M}_e$ (when $n$ is even) and $\mathcal{M}_o$ (when $n$ is odd) are invariant and asymptotically stable.

In the paper, because $\pi-W(\mathbf{\Gamma})$ is not differentiable in the whole state space,
we take the candidate Lyapunov function as
\[
V(\mathbf{\Gamma})
=2 \cos^2(W(\mathbf{\Gamma})/2)
= \max_{i\in\mathcal{V}} V_i(\mathbf{\Gamma})
\]
\revA{instead, where $V_i(\mathbf{\Gamma})= 2 \cos^2(\theta_{i,i+1}/2) = 1 + \Gamma_i^T \Gamma_{i+1}$.}
%The reason why $\pi-W(\mathbf{\Gamma})$ is not taken as a Lyapunov candidate is because %$\theta_{i,i+1}$ for any $i\in\mathcal{V}$
%it is not differentiable in the whole state space.
Calculating the time derivative of $V_i(\mathbf{\Gamma})$ along the trajectory of the closed-loop system (\ref{eq:closedloop}) yields
\begin{eqnarray}
\dot{V}_i(\mathbf{\Gamma})
&=& \Gamma_{i+1}^T \dot{\Gamma}_i + \Gamma_{i}^T \dot{\Gamma}_{i+1} \nonumber
\\
&=&
- \sin(\theta_{i,i+1}) k_{i,i+1}^T \sum_{j\in\mathcal{N}_i} \sin(\theta_{ij})k_{ij} \nonumber
\\
&&{}-\sin(\theta_{i+1,i}) k_{i+1,i}^T \sum_{j\in\mathcal{N}_{i+1}} \sin(\theta_{i+1,j})k_{i+1,j}.
%\\
%=&
%-  \sin(\theta_{i,i+1}) k_{i,i+1}^T \Bigg(\sum_{j\in\mathcal{N}_i}\sin(\theta_{ij})k_{ij}
%\nonumber\\[6pt]
%&\mspace{144mu} + \sum_{j\in\mathcal{N}_{i+1}} \sin(\theta_{j,i+1})k_{j,i+1}
\label{eq:dini:V}
%\Bigg)
\end{eqnarray}
We denote $D^{+} V(\mathbf{\Gamma})$ as the upper right Dini time derivative of $V(\mathbf{\Gamma})$ along the trajectory of the closed-loop system (\ref{eq:closedloop}).
(Details about Dini derivative can be found in \cite{giorgi,rouche}.)
And define the set $\mathcal{I}\subset\mathcal{V}$ of $\mathbf{\Gamma}\in(\mathcal{S}^2)^n$ as %$\mathcal{I}(\mathbf{\Gamma}) = \{ i\in\mathcal{V} : V(\mathbf{\Gamma}) = V_i(\mathbf{\Gamma})\}$, or equivalently,
\begin{eqnarray}
\mathcal{I}(\mathbf{\Gamma})
&=& \{ i\in\mathcal{V} : V(\mathbf{\Gamma}) = V_i(\mathbf{\Gamma}) \}
\nonumber \\
&=& \{ i\in\mathcal{V} : W(\mathbf{\Gamma}) = \theta_{i,i+1} \}.
\nonumber
\end{eqnarray}

\section{Antipodal formation}\label{sec:anfor}

In this section, we suppose the total number of agents is even and discuss the stability of the set $\mathcal{M}_e$ when $\mathcal{G}$ is a ring.
To this end, \revA{we first study a quantified relationship between the distance of the state $\mathbf{\Gamma}$ to the antipodal formation manifold $\mathcal{M}_e$ and the minimal geodesic distance of neighboring agents in Lemma~\ref{lem:even:bound}. This relationship contribute to the proof of Lyapunov stability of $\mathcal{M}_e$ in Theorem~\ref{thm:even:undirected:almost}.
 %At last, Theorem~\ref{thm:even:undirected:almost} gives the almost global stability in the undirected ring case.
 }

\begin{lem}\label{lem:even:bound}
Suppose $n$ is even and $\mathbf{\Gamma}\in(\mathcal{S}^2)^n$. Then
\begin{eqnarray}
\label{eq:lem:even:bound1}
&&W(\mathbf{\Gamma}) \geq \pi-2 \, d_{\mathcal{M}_e}(\mathbf{\Gamma}),
\\[6pt]
\label{eq:lem:even:bound2}
&&d_{\mathcal{M}_e}(\mathbf{\Gamma})
\leq \frac{n}{2} \left[\pi - W(\mathbf{\Gamma})\right].
\end{eqnarray}
\end{lem}

\begin{pf}
%The triangular inequality indicates that for any $\mathbf{\Gamma}\in (\mathcal{S}_2)^n$, and $i,j,k \in \mathcal{V}$, $d_{\mathcal{S}^2}(\Gamma_i, \Gamma_j) \leq d_{\mathcal{S}^2}(\Gamma_i, \Gamma_k)+d_{\mathcal{S}^2}(\Gamma_k, \Gamma_j)$.
Since there is a $\bar{\mathbf{\Gamma}}\in \mathcal{M}_e$ with $\bar{\Gamma}_i=-\bar{\Gamma}_{i+1}$ for any $i\in\mathcal{V}$
such that $\max_{i\in\mathcal{V}}d_{\mathcal{S}^2}(\Gamma_i, \bar{\Gamma}_i) = d_{\mathcal{M}_e}(\mathbf{\Gamma})$,
we obtain \revA{by the triangle inequality} that
\begin{eqnarray}
d_{\mathcal{S}^2}(\Gamma_i, \Gamma_{i+1})
&&\geq d_{\mathcal{S}^2}(\Gamma_i, \bar{\Gamma}_{i+1}) - d_{\mathcal{S}^2}(\bar{\Gamma}_{i+1}, \Gamma_{i+1})
\nonumber
\\
&&\geq d_{\mathcal{S}^2}(\bar{\Gamma}_i, \bar{\Gamma}_{i+1}) - d_{\mathcal{S}^2}(\bar{\Gamma}_{i}, \Gamma_{i}) - d_{\mathcal{M}_e}(\mathbf{\Gamma})
\nonumber \\
&&\geq \pi - 2\, d_{\mathcal{M}_e}(\mathbf{\Gamma}),
\quad \forall i\in\mathcal{V}, \nonumber
\end{eqnarray}
which implies (\ref{eq:lem:even:bound1}).

Let $\tilde{\Gamma}_i=(-1)^{i+\frac{n}{2}} \Gamma_{\frac{n}{2}}$ for any $i\in\mathcal{V}$.
Then $\tilde{\mathbf{\Gamma}}=\{\tilde{\Gamma}_i\}_{i\in\mathcal{V}} \in\mathcal{M}_e$.
\revA{
Since $d_{\mathcal{M}_e}(\mathbf{\Gamma})=\inf_{{\mathbf{\Gamma}}'\in\mathcal{M}_e }
d_{(\mathcal{S}^2)^n} (\mathbf{\Gamma},{\mathbf{\Gamma}}')
\leq d_{(\mathcal{S}^2)^n}(\mathbf{\Gamma},\tilde{\mathbf{\Gamma}})$,
we only need to show that $d_{(\mathcal{S}^2)^n}(\mathbf{\Gamma},\tilde{\mathbf{\Gamma}})
\leq \frac{n}{2} \left[\pi - W(\mathbf{\Gamma})\right].$
}
%Let $\tilde{\mathbf{\Gamma}}\in(\mathcal{S}^2)^n$ be the point such that $\tilde{\Gamma}_i=(-1)^{i+\frac{n}{2}} \Gamma_{\frac{n}{2}}$ for any $i\in\mathcal{V}$.

Notice that $d_{\mathcal{S}^2}(\Gamma_{\frac{n}{2}},\tilde{\Gamma}_{\frac{n}{2}})=0$.
Then, for any $i<n/2$,
\begin{eqnarray}
&&d_{\mathcal{S}^2} ((-1)^{i}\Gamma_i,(-1)^{\frac{n}{2}}\Gamma_{\frac{n}{2}})
\nonumber \\
\leq&&  \sum_{j=i}^{\frac{n}{2}-1} d_{\mathcal{S}^2} ((-1)^{j}\Gamma_j,(-1)^{j+1}\Gamma_{j+1})
\nonumber \\
=&& \sum_{j=i}^{\frac{n}{2}-1} \left[ \pi - d_{\mathcal{S}^2} (\Gamma_j,\Gamma_{j+1}) \right]
\leq \sum_{j=i}^{\frac{n}{2}-1} \left[ \pi - W(\mathbf{\Gamma}) \right].
\nonumber
\end{eqnarray}
Similarly, for any $i>n/2$,
%\begin{eqnarray}
%&\mspace{21mu}d_{\mathcal{S}^2} ((-1)^{\frac{n}{2}}\Gamma_{\frac{n}{2}}, (-1)^{i}\Gamma_i)
%\\
%&\leq  \sum_{j=\frac{n}{2}}^{i-1} d_{\mathcal{S}^2} ((-1)^{j}\Gamma_j,(-1)^{j+1}\Gamma_{j+1})
%\\
%&= \sum_{j=\frac{n}{2}}^{i-1} \left[ \pi - d_{\mathcal{S}^2} (\Gamma_j,\Gamma_{j+1}) \right]
%\leq \sum_{j=\frac{n}{2}}^{i-1} \left[ \pi - W(\mathbf{\Gamma}) \right].
%\end{eqnarray}
\[
d_{\mathcal{S}^2} ((-1)^{\frac{n}{2}}\Gamma_{\frac{n}{2}}, (-1)^{i}\Gamma_i)
\leq \sum_{j=\frac{n}{2}}^{i-1} \left[ \pi - W(\mathbf{\Gamma}) \right].
\]
Hence, $d_{(\mathcal{S}^2)^n}(\mathbf{\Gamma},\tilde{\mathbf{\Gamma}})
= \max_{i\in\mathcal{V}} d_{\mathcal{S}^2} (\Gamma_i,(-1)^{i+\frac{n}{2}}\Gamma_{\frac{n}{2}})
\leq \frac{n}{2} \left[\pi - W(\mathbf{\Gamma})\right]$.
% The proof is completed.
\qed
\end{pf}

\revA{
In the next lemma, the Dini time derivative of $V(\mathbf{\Gamma})$ is analyzed along the trajectory of the closed-loop system.
}

\begin{lem}\label{lem:even:diniv}
Suppose $n$ is even, Assumption \ref{assump:graph} holds and $W(\mathbf{\Gamma})\geq\pi-2\pi/n$. Then
$\dot{V}_i(\mathbf{\Gamma})\leq 0$ for any $i\in\mathcal{I}(\mathbf{\Gamma})$ and $D^{+} V(\mathbf{\Gamma}) \leq 0$.
Moreover, if $\dot{V}_i(\mathbf{\Gamma})=0$ for some $i\in\mathcal{I}(\mathbf{\Gamma})$, then
\begin{equation}\label{eq:anti:temp}
\sin(\theta_{i,i+1}) k_{i,i+1} = \sin(\theta_{i+1,i+2}) k_{i+1,i+2}.
\end{equation}
%and $D^{+} V(\mathbf{\Gamma})=0$ if and only if %$\mathbf{\Gamma}\in\mathcal{M}_e$ or
%\[
%\Gamma_i = \exp\left((i-1)d\frac{2\pi}{n}\widehat{u}\right) v,
%\quad \forall i\in\mathcal{V}
%\end{equation}
%for some $u,v\in\mathcal{S}^2$, $u^Tv=0$ and $d=\lfloor \frac{n}{4}\rfloor, \lfloor \frac{n}{4}\rfloor +1, \ldots, \lfloor \frac{n}{2}\rfloor$.
\end{lem}

\begin{pf}
Take any $\mathbf{\Gamma}$ satisfying $W(\mathbf{\Gamma})\geq \pi-2\pi/n$.

When $n=2$, $V(\mathbf{\Gamma})=V_1(\mathbf{\Gamma})=V_2(\mathbf{\Gamma})=2\cos^2(\theta_{12}/2)$.
By (\ref{eq:dini:V}), their time derivative along the closed-loop system (\ref{eq:closedloop}) is $-2\sin^2(\theta_{12})$ if $\mathcal{G}$ is undirected and $-\sin^2(\theta_{12})$ if $\mathcal{G}$ is directed, then the conclusions follow consequently.

When $n\geq 4$, take any $i\in\mathcal{I}(\mathbf{\Gamma})$, then for any $j\in\mathcal{V}$,
$\pi/2 \leq \theta_{i,i+1} \leq \theta_{j,j+1}$ , which implies
\begin{equation}\label{eq:2016092905}
\sin(\theta_{i,i+1}) \geq \sin(\theta_{j,j+1}),
\end{equation}
for any $j\in\mathcal{V}$.
If $\mathcal{G}$ is an undirected ring, \revA{since $k_{i,j}^Tk_{i,k} \in [-1,1]$ for any $i,j,k \in \mathcal{V}$, we have}
\begin{eqnarray}
\dot{V}_i(\mathbf{\Gamma})
&=&
- \sin(\theta_{i,i+1})  \big[2\sin(\theta_{i,i+1})  \nonumber \\
&&{} +\sin(\theta_{i,i-1}) k_{i,i+1}^Tk_{i,i-1} \nonumber\\
&&{}+\sin(\theta_{i+1,i+2}) k_{i+1,i}^T k_{i+1,i+2} \big]
\leq 0. \nonumber
\end{eqnarray}
Then, if $\mathcal{G}$ is a directed ring,
\begin{eqnarray}
\dot{V}_i(\mathbf{\Gamma})
&=&
- \sin(\theta_{i,i+1})  \big[ \sin(\theta_{i,i+1}) \nonumber\\
&&{}+  \sin(\theta_{i+1,i+2}) k_{i+1,i}^T k_{i+1,i+2} \big] \leq 0. \nonumber
\end{eqnarray}
In both cases,
$D^{+} V(\mathbf{\Gamma}) = \max_{i\in\mathcal{I}(\mathbf{\Gamma})} \dot{V}_i(\mathbf{\Gamma}) \leq 0$.
\revA{Moreover, since (\ref{eq:2016092905}) is fulfilled for any $j\in\mathcal{V}$, $\dot{V}_i(\mathbf{\Gamma})=0$ implies (\ref{eq:anti:temp}).}
\qed
\end{pf}

Due to \rev{the fact}  $\mathcal{M}_e$ is an equilibrium set of the closed-loop system (\ref{eq:closedloop}),
it is (positively) invariant.
We are now ready to give the local stability of the invariant set $\mathcal{M}_e$.

\begin{thm}\label{thm:even}
Suppose $n$ is even. Then under Assumption~\ref{assump:graph}, the
invariant set $\mathcal{M}_e$ is asymptotically stable with the region
of attraction containing
%\[
%\Omega_e = \left\{ \mathbf{\Gamma}\in(\mathcal{S}^2)^n :
%d_{\mathcal{M}_e}(\mathbf{\Gamma}) < \pi/n \right\}.
%\end{equation}
\[
\Omega_e = \left\{ \mathbf{\Gamma}\in(\mathcal{S}^2)^n :
W(\mathbf{\Gamma})> \pi-2\pi/n \right\}.
\]
\end{thm}

\begin{pf}
\revA{In order to prove the asymptotical stability of the invariant set $\mathcal{M}_e$, both stability and atrractivity of $\mathcal{M}_e$ need to be shown.}
We first show that $\mathcal{M}_e$ is stable.
For any $\varepsilon\in(0,\pi]$, let $\delta=\varepsilon/n$.
Then for any initial state $\mathbf{\Gamma}(0)$ satisfies
$d_{\mathcal{M}_e}(\mathbf{\Gamma}(0))<\delta$,
$W(\mathbf{\Gamma}(0))>\pi-2\delta$ by (\ref{eq:lem:even:bound1}).
Due to \rev{the fact} $\pi-2\delta\geq \pi-2\pi/n$, $D^{+} V(\mathbf{\Gamma}(0)) \leq 0$ by Lemma~\ref{lem:even:diniv}.
Then $\mathbf{\Gamma}(t)\in\Omega_e$ for all $t\geq 0$ and $V(\mathbf{\Gamma}(\cdot))$ is non-increasing on $[0,\infty)$.
%According to the comparison lemma (see Lemma~3.4 in \cite{khalil}), we obtain $V(\mathbf{\Gamma}(t)) \leq V(\mathbf{\Gamma}(0))$ for any $t\geq 0$.
%According to the comparison lemma (see Lemma~3.4 in \cite{khalil}), we obtain $V(\mathbf{\Gamma}(t)) \leq V(\mathbf{\Gamma}(0))$ for any $t\geq 0$.
Therefore,
$W(\mathbf{\Gamma}(t)) \geq W(\mathbf{\Gamma}(0))>\pi-2\delta$ for any $t\geq 0$,
which implies $d_{\mathcal{M}_e}(\mathbf{\Gamma}(t))
%\leq \frac{n}{2} \left[ \pi - W(\mathbf{\Gamma}(t)) \right]
<n \delta = \varepsilon$ for any $t\geq 0$ by (\ref{eq:lem:even:bound2}).
Hence, $\mathcal{M}_e$ is stable.
%Meanwhile, this also shows the set $\Omega_e$ is positively invariant.

In the following, we use $\mathbf{\Gamma}(t;\bar{\mathbf{\Gamma}})$ to denote the trajectory of the closed-loop system (\ref{eq:closedloop}) starting from $\bar{\mathbf{\Gamma}} \in (\mathcal{S}^2)^n$ at $t=0$.

\revA{Next we prove that $\mathcal{M}_e$ is also attractive. To this end, we only need to show, for any $\bar{\mathbf{\Gamma}}\in\Omega_e$,
\begin{equation}\label{eq:even:temp1}
\lim_{t\rightarrow\infty} d_{\mathcal{M}_e}(\mathbf{\Gamma}(t;\bar{\mathbf{\Gamma}})) = 0.
\end{equation}
}

First, from the above proof, we have $D^{+}V(\mathbf{\Gamma}(t;\bar{\mathbf{\Gamma}}))\leq 0$ for any $t\geq 0$. In addition, since $\bar{\mathbf{\Gamma}}\in\Omega_e$,
\[
V(\mathbf{\Gamma}(0;\bar{\mathbf{\Gamma}})) =2\cos^2(W(\bar{\mathbf{\Gamma}})/2)
<2\cos^2((\pi-2\pi/n)/2).
\]

Combining with $V(\mathbf{\Gamma}(t;\bar{\mathbf{\Gamma}}))\geq 0$ and
$D^{+} V(\mathbf{\Gamma}(t;\bar{\mathbf{\Gamma}})) \leq 0$ for any $t\geq 0$,
we infer there is a constant
\[
V^*=2\cos^2(\alpha^*/2), \quad\mathrm{where}\;\alpha^*\in(\pi-2\pi/n,\pi],
\]
such that $\lim_{t\rightarrow\infty} V(\mathbf{\Gamma}(t;\bar{\mathbf{\Gamma}}))=V^*$.
If $\alpha^*=\pi$, then (\ref{eq:even:temp1}) follows consequently by Lemma~\ref{lem:basic:even}.
In the following, we suppose $\alpha^*\in(\pi-2\pi/n,\pi)$ and show this leads to a contradiction.

Since the trajectory $\mathbf{\Gamma}(t;\bar{\mathbf{\Gamma}})$ is bounded, its positive limit set $L(\bar{\mathbf{\Gamma}})\subset\Omega_e$ is nonempty and positively invariant.
For any limit point $\tilde{\mathbf{\Gamma}}\in L(\bar{\mathbf{\Gamma}})$, there exists an increasing time sequence $\{t_k\}$ such that as $k\rightarrow\infty$, $t_k\rightarrow\infty$ and $\mathbf{\Gamma}(t_k;\bar{\mathbf{\Gamma}})\rightarrow\tilde{\mathbf{\Gamma}}$.
By continuity of $V$, we obtain
\begin{equation}\label{eq:even:temp2}
V(\tilde{\mathbf{\Gamma}})
% = V(\lim_{k\rightarrow\infty}\mathbf{\Gamma}(t_k;\bar{\mathbf{\Gamma}}))
= V^*,\quad
W(\tilde{\mathbf{\Gamma}}) = \alpha^*,
\quad \forall \tilde{\mathbf{\Gamma}}\in L(\bar{\mathbf{\Gamma}}).
\end{equation}
Take any $\tilde{\mathbf{\Gamma}}\in L(\bar{\mathbf{\Gamma}})$ and consider another trajectory $\mathbf{\Gamma}(t;\tilde{\mathbf{\Gamma}})$. %starting from $\tilde{\mathbf{\Gamma}}$.
Clearly, $D^{+}V(\mathbf{\Gamma}(t;\tilde{\mathbf{\Gamma}}))\leq 0$ and $V(\mathbf{\Gamma}(t;\tilde{\mathbf{\Gamma}}))\leq V^{*}$ for any $t\geq 0$.
Define
\[
\mathcal{K}(t) %= \{ i\in\mathcal{V}: V_i(\mathbf{\Gamma}(t;\tilde{\mathbf{\Gamma}}))=V^* \}
= \{ i\in\mathcal{V}: \theta_{i,i+1}(t;\tilde{\mathbf{\Gamma}}) = \alpha^*
\}.
\]
Since $W(\tilde{\mathbf{\Gamma}}) = \alpha^*$, $\mathcal{K}(0)\neq \emptyset$.
Then we claim that $\mathcal{K}(t)=\emptyset$ for any $t>0$.
Suppose this is not true, then there exist $i\in\mathcal{V}$ and $s>0$ such that $i\in\mathcal{K}(s)$.
Hence, $i\in\mathcal{I}(\mathbf{\Gamma}(s;\tilde{\mathbf{\Gamma}}))$, and by Lemma~\ref{lem:even:diniv} there are two cases for $\dot{V}_i(\mathbf{\Gamma}(s;\tilde{\mathbf{\Gamma}}))$ which are discussed as follows.

Case 1: $\dot{V}_i(\mathbf{\Gamma}(s;\tilde{\mathbf{\Gamma}}))<0$.
Since $\dot{V}_i(\mathbf{\Gamma}(t;\tilde{\mathbf{\Gamma}}))$ is continuous in time,
there exists an $s'\in(0,s)$ such that $\dot{V}_i(\mathbf{\Gamma}(t;\tilde{\mathbf{\Gamma}})) <0$ for any $t\in(s',s]$.
Combining with $V_i(\mathbf{\Gamma}(s';\tilde{\mathbf{\Gamma}})) \leq V^*$, we obtain $V_i(\mathbf{\Gamma}(s;\tilde{\mathbf{\Gamma}})) < V^*$, which contradicts $i\in\mathcal{K}(s)$.

Case 2: $\dot{V}_i(\mathbf{\Gamma}(s;\tilde{\mathbf{\Gamma}}))=0$.
Let $u=k_{i,i+1}(s;\tilde{\mathbf{\Gamma}})$.
By Lemma~\ref{lem:even:diniv}, since $\alpha^*\in(\pi-2\pi/n,\pi)$, we obtain
\[
\theta_{i+1,i+2}(s;\tilde{\mathbf{\Gamma}}) = \alpha^*,
\quad k_{i+1,i+2}(s;\tilde{\mathbf{\Gamma}}) = u.
\]
Therefore, $i+1\in\mathcal{K}(s)$.
If $\dot{V}_{i+1}(\mathbf{\Gamma}(s;\tilde{\mathbf{\Gamma}}))<0$, this turns to Case 1 and leads to a contradiction.
Thus, $\dot{V}_{i+1}(\mathbf{\Gamma}(s;\tilde{\mathbf{\Gamma}}))=0$.
Repeating these arguments, since $\mathcal{G}$ is a ring, we obtain
\[
\theta_{i,i+1}(s;\tilde{\mathbf{\Gamma}})=\alpha^*, \quad
k_{i,i+1}(s;\tilde{\mathbf{\Gamma}}) = u, \quad
\forall i\in\mathcal{V}.
\]
This indicates that $\Gamma_1(s;\tilde{\mathbf{\Gamma}}),\ldots,\Gamma_n(s;\tilde{\mathbf{\Gamma}})$ lie on a great circle whose axis is $u$, and $\Gamma_{i+1}(s;\tilde{\mathbf{\Gamma}})=\exp(\alpha^*\widehat{u})\Gamma_i(s;\tilde{\mathbf{\Gamma}})$ for any $i\in\mathcal{V}$.
Hence, $\Gamma_1(s;\tilde{\mathbf{\Gamma}})
= \exp(n\alpha^*\widehat{u})\Gamma_1(s;\tilde{\mathbf{\Gamma}})$,
which is a contradiction since rotating
$\Gamma_1(s;\tilde{\mathbf{\Gamma}})$ about $u$ that is orthogonal to
$\Gamma_1(s;\tilde{\mathbf{\Gamma}})$ through an angle
$n\alpha^*\in((n-2)\pi,n\pi)$ can not return to
$\Gamma_1(s;\tilde{\mathbf{\Gamma}})$.

As a result, we conclude that $\mathcal{K}(t)=\emptyset$ for any $t>0$.
Then $V(\mathbf{\Gamma}(t;\tilde{\mathbf{\Gamma}}))<V^*$ for any $t>0$.
This contradicts (\ref{eq:even:temp2}) since $L(\bar{\mathbf{\Gamma}})$ is positively invariant.
Hence, $\alpha^*=\pi$ and the conclusion follows consequently. \qed
\end{pf}

Since $\mathcal{M}_e$ is an equilibrium set of the closed-loop system (\ref{eq:closedloop}),
the above theorem ensures that the reduced attitudes of the agents with initial values in $\Omega_e$ reach a static antipodal formation eventually provided that the inter-agent graph is a ring, either undirected or directed.
When $n=2$, $\Omega_e= \{ \mathbf{\Gamma}\in (\mathcal{S}^2)^2 : \Gamma_1 \neq \Gamma_2\}$ has measure one and is everywhere dense, which means the domain of attraction of $\mathcal{M}_e= \{ \mathbf{\Gamma}\in (\mathcal{S}^2)^2 : \Gamma_1 =-\Gamma_2\}$ is almost all $(\mathcal{S}^2)^2$.
%=(\mathcal{S}^2)^2 \backslash

\begin{cor}\label{corollary}
For the trivial case $n=2$, $\mathcal{M}_e$ is almost globally asymptotically stable under Assumption~\ref{assump:graph}.
\end{cor}

Because the state space $(\mathcal{S}^2)^n$ is compact without boundary, there are multiple disjoint equilibrium manifolds for the closed-loop system~(\ref{eq:closedloop}).
The next theorem shows that the region of attraction for the antipodal formation manifold $\mathcal{M}_e$ is almost all $(\mathcal{S}^2)^n$ provided that $\mathcal{G}$ is an undirected ring.

\begin{thm}\label{thm:even:undirected:almost}
Suppose $n$ is even and $\mathcal{G}$ is an undirected ring.
As time tends to infinity, every trajectory of the closed-loop system $(\ref{eq:closedloop})$ approaches one of its equilibrium manifolds,
among which $\mathcal{M}_e$ is the only stable one.
\end{thm}

\begin{pf}
\revA{
By Corollary \ref{corollary}, we only need to consider the nontrivial case $n>2$ in the following proof.}

\revA{
First, we show every trajectory of (\ref{eq:closedloop}) approaches the equilibrium set of the closed-loop system.
}Let $\bar{V}(\mathbf{\Gamma})=\sum_{i=1}^n V_i(\mathbf{\Gamma}) $.
By (\ref{eq:dini:V}), the time derivative of $\bar{V}(\mathbf{\Gamma})$ along the trajectory of the closed-loop system~(\ref{eq:closedloop}) is given as
\[
\dot{\bar{V}}(\mathbf{\Gamma})
=
- \sum_{i=1}^n \left\| \sum_{j\in\mathcal{N}_i} \sin(\theta_{ij})k_{ij} \right\|^2.
\]
Thus, by LaSalle's invariance principle \cite{khalil}, as $t\rightarrow\infty$ every trajectory of (\ref{eq:closedloop}) approaches
$\mathcal{D} = \{ \mathbf{\Gamma}\in(\mathcal{S}^2)^n :
 \sin(\theta_{i,i-1})k_{i-1,i} = \sin(\theta_{i,i+1})k_{i,i+1}, \forall i\in\mathcal{V} \}$,
which is the entire equilibrium set of the closed-loop system.
By Theorem~\ref{thm:even}, the equilibrium set $\mathcal{M}_e\subset\mathcal{D}$ is asymptotically stable.

\revA{
Next we take any $\mathbf{\Gamma}^*\in\mathcal{D}\backslash\mathcal{M}_e$ and show that $\mathbf{\Gamma}^*$ is unstable by applying the Lyapunov's indirect method.}
Denote $\{(\psi_i^*,\varphi_i^*)\}_{i\in\mathcal{V}}$ as the parametrization of $\mathbf{\Gamma}^*$ corresponding to (\ref{eq:gamma:para}).
Due to \rev{the fact} $\mathbf{\Gamma}^*\in\mathcal{D}$, there is a $u\in\mathcal{S}^2$ such that $\Gamma_1^*,\ldots,\Gamma_n^*$ lie on the great circle of $\mathcal{S}^2$ orthogonal to $u$.
Denote $e_3=[0,0,1]^T\in\mathcal{S}^2$ as the north pole.
If $u\neq \pm e_3$, let $\alpha=\arccos(u^Te_3)$ and $v=\widehat{u}e_3/\sin(\alpha)$.
Then through coordinate transformation $\bar{\Gamma}_i=\exp(\alpha\widehat{v})\Gamma_i$ for any $i\in\mathcal{V}$, the corresponding equilibrium $\bar{\Gamma}_1^*,\ldots,\bar{\Gamma}_n^*$ of the new system lie on the great circle of $\mathcal{S}^2$ orthogonal to $e_3$, i.e., the equator.
And by Lemma~\ref{lem:invariant}, the stability of $\mathbf{\Gamma}^*$ with respect to the original system is consistent with the stability of $\mathbf{\bar{\Gamma}}^*$ with respect to the new system.
Therefore, without loss of generality, we assume $u=e_3$.
Then $\varphi_i^*=0$ and $\cos(\psi_i^*-\psi_j^*)=\cos(\theta_{ij}^*)$ for any $i,j\in\mathcal{V}$.

Let $\bar{\psi}_i = \psi_i-\psi_i^*$.
Computing the linearization of the closed-loop system (\ref{eq:closedloop:psi})-(\ref{eq:closedloop:varphi}) about $\{(\psi_i^*,\varphi_i^*)\}_{i\in\mathcal{V}}$ yields
\begin{eqnarray}
\dot{\bar\psi}_i &=& \sum_{j\in\mathcal{N}_i} \cos(\theta_{ij}^*)  (\bar{\psi}_i - \bar{\psi}_j),
\quad i=1,2,\ldots,n,
\label{eq:closedloop:psi:linear}\\[3pt]
\dot{\varphi}_i
&=& \sum_{j\in\mathcal{N}_i} \left[ \cos(\theta_{ij}^*) \varphi_i - \varphi_j \right],
\quad i=1,2,\ldots,n.
\label{eq:closedloop:varphi:linear}
\end{eqnarray}
The above linearizations for $\{\psi_i\}_{i\in\mathcal{V}}$ and $\{\varphi_i\}_{i\in\mathcal{V}}$ are decoupled and we then denote the respective system matrices of (\ref{eq:closedloop:psi:linear}) and (\ref{eq:closedloop:varphi:linear}) as $A_{\psi}(\mathbf{\Gamma}^*), A_{\varphi}(\mathbf{\Gamma}^*) \in\mathbb{R}^{n\times n}$,
both of which are symmetric due to \rev{the fact} $\mathcal{G}$ is undirected.
Let $x_i=(-1)^i$ for any $i\in\mathcal{V}$ and $\boldsymbol{x}=[x_1,\ldots,x_n]^T$.
Then
\[
\boldsymbol{x}^T A_{\varphi}(\mathbf{\Gamma}^*) \boldsymbol{x}
= \sum_{i=1}^n \left[ 2 + \cos(\theta_{i,i+1}^*) + \cos(\theta_{i,i-1}^*) \right].
\]
Due to \rev{the fact} $\mathbf{\Gamma}^*\in\mathcal{D}\backslash\mathcal{M}_e$, there is an $i\in\mathcal{V}$ such that $\theta_{i,i+1}^*<\pi$, which implies $\boldsymbol{x}^T A_{\varphi}(\mathbf{\Gamma}^*) \boldsymbol{x}>0$ and $A_{\varphi}(\mathbf{\Gamma}^*)$ has at least one positive eigenvalue, i.e., $\mathbf{\Gamma}^*$ is unstable.
\qed
\end{pf}

\begin{rem}\label{rem:even:eigenvalue}
When $\mathcal{G}$ is a directed ring, all trajectories of the closed-loop system (\ref{eq:closedloop}) eventually approach an invariant set which also contains non-equilibrium points.
When $\mathcal{G}$ is an undirected ring, the Jacobian matrix of the closed-loop system (\ref{eq:closedloop:psi})-(\ref{eq:closedloop:varphi}) about any equilibrium $\mathbf{\Gamma}^*\in \mathcal{M}_e$ has $2$ zero eigenvalues and $2n-2$ negative eigenvalues.
The existence of the zero eigenvalues makes the linearization fail to determine the stability of  $\mathbf{\Gamma}^*$.
As a fact, any eigenvector of the zero eigenvalue corresponds to a tangent vector of $\mathcal{M}_e$ at $\mathbf{\Gamma}^*$.
An explanation for why there are $2$ zero eigenvalues is that the invariant manifold $\mathcal{M}_e$ has two degrees of freedom in the whole state space $(\mathcal{S}^2)^n$.
\end{rem}

\section{Cyclic formation}\label{sec:odd:cyclic}

In this section, we suppose the total number of agents is odd and discuss the stability of the set $\mathcal{M}_o$ when $\mathcal{G}$ is a ring.
%To proceed this subsection,
This section will proceed similar to Section~\ref{sec:anfor}, however, the theoretical analysis is more complicated.

\revA{
We start with a lemma to provide the upper bound of $d_{\mathcal{S}^2} (\Gamma_i,\bar{\Gamma}_{i})$, whose proof can be found in Appendix~\ref{appen:sec:proof}.
}

\begin{lem}\label{lem:odd:bound:0}
Suppose $n$ is odd, $\mathbf{\Gamma}\in(\mathcal{S}^2)^n$ and
\begin{equation}\label{eq:odd:temp1}
\left| \pi-\pi/n -\theta_{i,i+1} \right| \leq \nu^2,
\quad \forall i\in\mathcal{V}.
\end{equation}
If $\nu\in[0,\sqrt{2}/n]$, then $d_{\mathcal{M}_o}(\mathbf{\Gamma})\leq 2n\nu$.
\end{lem}

\revA{Next, in order to prove the Lyapunov stability of invariant set $\mathcal{M}_o$, the relationship between $W(\mathbf{\Gamma})$ and $d_{\mathcal{M}_o}(\mathbf{\Gamma})$ is investigated in Lemma~\ref{lem:odd:bound}.}

\begin{lem}\label{lem:odd:bound}
Suppose $n$ is odd and $\mathbf{\Gamma}\in(\mathcal{S}^2)^n$. Then
\begin{equation}\label{eq:lem:odd:bound1}
W(\mathbf{\Gamma}) \geq \pi-\pi/n-2\, d_{\mathcal{M}_o}(\mathbf{\Gamma}).
\end{equation}
If $\sqrt{4n^3(\pi-\pi/n-W(\mathbf{\Gamma}))} \leq 2\sqrt{2}$, then
\begin{equation}\label{eq:lem:odd:bound2}
d_{\mathcal{M}_o}(\mathbf{\Gamma}) \leq \sqrt{4n^3(\pi-\pi/n-W(\mathbf{\Gamma}))} .
\end{equation}
\end{lem}

\begin{pf}
The proof for (\ref{eq:lem:odd:bound1}) is similar to that of Lemma~\ref{lem:even:bound} and is omitted here.
%Denote $\xi = d_{\mathcal{M}_o}(\mathbf{\Gamma})$.
%Then there exists a $\bar{\mathbf{\Gamma}}\in \mathcal{M}_o$
%with $d(\bar{\Gamma}_i,\bar{\Gamma}_{i+1})=\pi-\pi/n$, $\forall i\in\mathcal{V}$
%such that $d_{\mathcal{S}^2}(\Gamma_i, \bar{\Gamma}_i) \leq \xi$, $\forall i\in\mathcal{V}$.
%Hence,
%\begin{eqnarray}
%d_{\mathcal{S}^2}(\Gamma_i, \Gamma_{i+1})
%&\geq& d_{\mathcal{S}^2}(\Gamma_i, \bar{\Gamma}_{i+1}) - d_{\mathcal{S}^2}(\bar{\Gamma}_{i+1}, \Gamma_{i+1})
%\nonumber \\
%&\geq& d_{\mathcal{S}^2}(\bar{\Gamma}_i, \bar{\Gamma}_{i+1}) - d_{\mathcal{S}^2}(\bar{\Gamma}_{i}, \Gamma_{i}) - \xi
%\nonumber \\
%&\geq& \pi-\pi/n - 2\xi,
%\quad \forall i\in\mathcal{V}. \nonumber
%\end{eqnarray}
%This implies $W(\mathbf{\Gamma})\geq \pi-\pi/n-2\xi$.

Let $\kappa = \sqrt{4n(\pi-\pi/n-W(\mathbf{\Gamma}))}$. Then
\[
\theta_{i,i+1} \geq W(\mathbf{\Gamma}) = \pi-\pi/n -\kappa^2/(4n),
\quad \forall i\in\mathcal{V}.
\]
Since
$(n-1)\pi \geq \sum_{i\in\mathcal{V}}\theta_{i,i+1}
\geq \max_{i\in\mathcal{V}} \theta_{i,i+1} + (n-1)W(\mathbf{\Gamma})$,
we obtain
\[
\max_{i\in\mathcal{V}} \theta_{i,i+1} \leq \pi-\pi/n + (n-1)\kappa^2/(4n).
\]
Therefore,
$| \pi - \pi/n - \theta_{i,i+1} | \leq \left(\kappa/2\right)^2$ for any $i\in\mathcal{V}$.
Since $\kappa/2 \in[0,\sqrt{2}/n]$,
$d_{\mathcal{M}_o}(\mathbf{\Gamma}) \leq n\kappa$ by Lemma~\ref{lem:odd:bound:0}.
\qed
\end{pf}

\begin{lem}\label{lem:odd:diniv}
Suppose $n$ is odd, Assumption \ref{assump:graph} holds and
$W(\mathbf{\Gamma})\geq \max\{\pi-3\pi/n, \pi/2 \}$.
Then
$\dot{V}_i(\mathbf{\Gamma})\leq 0$ for any $i\in\mathcal{I}(\mathbf{\Gamma})$ and
$D^{+} V(\mathbf{\Gamma}) \leq 0$.
If $\dot{V}_i(\mathbf{\Gamma})=0$ for some $i\in\mathcal{I}(\mathbf{\Gamma})$, then
$\sin(\theta_{i,i+1}) k_{i,i+1} = \sin(\theta_{i+1,i+2}) k_{i+1,i+2}$.
\end{lem}

\begin{pf}
For any $\mathbf{\Gamma}\in(\mathcal{S}^2)^n$ such that $W(\mathbf{\Gamma})\geq \max\{\pi-3\pi/n, \pi/2 \}$, it holds that
$\pi/2 \leq \theta_{i,i+1} \leq \theta_{j,j+1}$ for any $i\in\mathcal{I}(\mathbf{\Gamma})$ and $j\in\mathcal{V}$.
The remaining of the proof is similar to that of Lemma~\ref{lem:even:diniv} and is omitted here.\qed
\end{pf}

When Assumption~\ref{assump:graph} holds,
$\dot{V}_i(\mathbf{\Gamma}) = 0$ for any $i\in\mathcal{V}$ and
$\mathbf{\Gamma}\in\mathcal{M}_o$ by (\ref{eq:dini:V}).
Therefore, the set $\mathcal{M}_o$ is (positively) invariant with respect to the closed-loop system (\ref{eq:closedloop}).
In particular, $\mathcal{M}_o$ is an equilibrium set when the inter-agent graph $\mathcal{G}$ is an undirected ring.
With Lemma~\ref{lem:odd:bound} and Lemma~\ref{lem:odd:diniv},
we are ready to state the local stability of the invariant set $\mathcal{M}_o$ in the next theorem.

\begin{thm}\label{thm:odd}
Suppose $n$ is odd and Assumption~\ref{assump:graph} holds.
The invariant set $\mathcal{M}_o$ is asymptotically stable with the region
of attraction containing
\[
\Omega_o = \left\{ \mathbf{\Gamma}\in(\mathcal{S}^2)^n :
W(\mathbf{\Gamma})> \max\{\pi-3\pi/n, \pi/2 \} \right\}.
\]
\end{thm}

\begin{pf}
We first show that $\mathcal{M}_o$ is stable.
For any $\varepsilon\in(0,\pi]$, let $\delta=\min\{ \varepsilon^2/(8n^3),  1/n^3\}$.
Take any initial state $\mathbf{\Gamma}(0)$ satisfying
$d_{\mathcal{M}_o}(\mathbf{\Gamma}(0))<\delta$.
Then $W(\mathbf{\Gamma}(0))> \pi-\pi/n-2\delta$ by (\ref{eq:lem:odd:bound1}).
Since $\pi-\pi/n-2\delta
> \max\{\pi-3\pi/n, \pi/2 \}$,
$D^{+} V(\mathbf{\Gamma}(0)) \leq 0$ by Lemma~\ref{lem:odd:diniv},
implying $\mathbf{\Gamma}(t)\in\Omega_o$ for any $t\geq 0$ and $V(\mathbf{\Gamma}(\cdot))$ is non-increasing on $[0,\infty)$.
Hence,
$W(\mathbf{\Gamma}(t)) \geq W(\mathbf{\Gamma}(0))>\pi-\pi/n-2\delta$ for any $t\geq 0$.
Then due to \rev{the fact}
$\sqrt{4n^3(\pi-\pi/n-W(\mathbf{\Gamma}(t)))}
< \sqrt{8n^3 \delta}
\leq 2\sqrt{2}$ for any $t\geq 0$,
$d_{\mathcal{M}_o}(\mathbf{\Gamma}(t))
< \sqrt{8n^3 \delta} \leq \varepsilon$ for any $t\geq 0$ by (\ref{eq:lem:odd:bound2}).
Therefore, $\mathcal{M}_o$ is stable.
%Meanwhile, this also shows the set $\Omega_o$ is positively invariant.

%In the following, we use $\mathbf{\Gamma}(t;\bar{\mathbf{\Gamma}})$ to denote the trajectory of the closed-loop system (\ref{eq:closedloop}) starting from $\bar{\mathbf{\Gamma}} \in (\mathcal{S}^2)^n$.

Next we show that $\mathcal{M}_o$ is also attractive.
Denote $\mathbf{\Gamma}(t;\bar{\mathbf{\Gamma}})$ as the trajectory of the closed-loop system (\ref{eq:closedloop}) starting from $\bar{\mathbf{\Gamma}} \in (\mathcal{S}^2)^n$ at $t=0$.
Take any $\bar{\mathbf{\Gamma}}\in\Omega_o$.
By the above proof, $D^{+}V(\mathbf{\Gamma}(t;\bar{\mathbf{\Gamma}}))\leq 0$ for any $t\geq 0$.
Then since
\[
V(\mathbf{\Gamma}(0;\bar{\mathbf{\Gamma}})) =2\cos^2(W(\bar{\mathbf{\Gamma}})/2)
<2\cos^2((\pi-3\pi/n)/2)
\]
and
$V(\mathbf{\Gamma}(t;\bar{\mathbf{\Gamma}}))\geq 2\cos^2((\pi-\pi/n)/2)$ for any $t\geq 0$ by Lemma~\ref{lem:basic:odd}, there is a constant
\[
V^*=2\cos^2(\alpha^*/2),
\quad\mathrm{where}\;
\alpha^*\in(\pi-3\pi/n,\pi-\pi/n],
\]
such that $\lim_{t\rightarrow\infty} V(\mathbf{\Gamma}(t;\bar{\mathbf{\Gamma}}))=V^*$.
Similar to the proof of Theorem \ref{thm:even}, it can be verified that $\alpha^{*}=\pi-\pi/n$, which implies
$\lim_{t\rightarrow\infty} d_{\mathcal{M}_o}(\mathbf{\Gamma}(t;\bar{\mathbf{\Gamma}})) = 0$
by Lemma~\ref{lem:basic:odd}. The proof is completed.
\qed
\end{pf}

When the inter-agent graph $\mathcal{G}$ is an undirected ring, since
$\mathcal{M}_o$ is an equilibrium set of the closed-loop system
(\ref{eq:closedloop}), the reduced attitudes of the system with initial values in the set $\Omega_o$ reach a static cyclic formation eventually by Theorem~\ref{thm:odd}.
\rev{
When $\mathcal{G}$ is a directed ring, a rotating cyclic reduced
attitude formation is attained, where the whole formation is rotating
about the axis of the great circle. Moreover, by the control law (\ref{eq:controller}), the magnitude of the angular velocity of any agent $i\in\mathcal{V}$ equals $\sin(\theta_{i,i+1}) = \sin(\pi-\pi/n) = \sin(\pi/n)$.
}

The next theorem further demonstrates that the region of attraction for the cyclic formation manifold $\mathcal{M}_o$ is almost all $(\mathcal{S}^2)^n$ provided that $\mathcal{G}$ is an undirected ring.

\begin{thm}\label{thm:odd:undirected:almost}
Suppose $n$ is odd and $\mathcal{G}$ is an undirected ring.
As time tends to infinity, every trajectory of the closed-loop system $(\ref{eq:closedloop})$ approaches one of its equilibrium manifolds,
among which $\mathcal{M}_o$ is the only stable one.
\end{thm}

\begin{pf}
Let $\mathcal{D} = \{ \mathbf{\Gamma}\in(\mathcal{S}^2)^n :
 \sin(\theta_{i,i-1})k_{i-1,i} = \sin(\theta_{i,i+1})k_{i,i+1}, \forall i\in\mathcal{V} \}$ be the entire equilibrium set of the closed-loop system~(\ref{eq:closedloop})
and take any $\mathbf{\Gamma}^*\in\mathcal{D}\backslash\mathcal{M}_o$.
Similar to the proof of Theorem~\ref{thm:even:undirected:almost},
we only need to show that $\mathbf{\Gamma}^*$ is unstable.
To this end, we analyze the eigenvalues of the Jacobian matrix $A_{\psi}(\mathbf{\Gamma}^*)$ and $A_{\varphi}(\mathbf{\Gamma}^*)$ of the linearization system (\ref{eq:closedloop:psi:linear})-(\ref{eq:closedloop:varphi:linear}).

Denote $\{(\psi_i^*,\varphi_i^*)\}_{i\in\mathcal{V}}$ as the parametrization of $\mathbf{\Gamma}^*$ corresponding to (\ref{eq:gamma:para}).
Without loss of generality, we assume $\varphi_i^*=0$ and $\cos(\psi_i^*-\psi_j^*)=\cos(\theta_{ij}^*)$ for any $i,j\in\mathcal{V}$.
Due to \rev{the fact} $\mathbf{\Gamma}^*\in\mathcal{D}$, there exist $a\in[0,1]$ and $u\in\mathcal{S}^2$ such that $\sin(\theta_{i,i+1}^*)=a$ and $k_{i,i+1}^*=u$ for any $i\in\mathcal{V}$. We group the equilibrium $\mathbf{\Gamma}^*$ into two categories and discuss them as follows.

%Case 1: $a=0$, i.e., $\theta_{i,i+1}^*$ is either $0$ or $\pi$ for any $i\in\mathcal{V}$.
%By Lemma~\ref{lem:basic:odd}, there exists at one node $i\in\mathcal{V}$ such that $\theta_{i,i+1}^*=0$.
%Let $\boldsymbol{x}=[x_1,\ldots,x_n]^T\in\mathbb{R}^n$ such that $x_i=1$, $x_{i+1}=-1$ and other entries are $0$.
%Then
%$\boldsymbol{x}^T A_{\psi}(\mathbf{\Gamma}^*) \boldsymbol{x}
%= 4 + \cos(\theta_{i,i-1}^*) + \cos(\theta_{i+1,i+2}^*)>0$.

Category 1: $\theta_{12}^*=\theta_{23}^*=\cdots=\theta_{n1}^*$ does not hold,
that is, $a\in[0,1)$ and  $\theta_{i-1,i}^*+\theta_{i,i+1}^*=\pi$ for some $i\in\mathcal{V}$.
If $\theta_{i-1,i}^*<\pi/2$, let $\boldsymbol{x}=[x_1,\ldots,x_n]^T\in\mathbb{R}^n$ such that $x_i=-x_{i-1}=1$ and other entries be $0$,
then
$\boldsymbol{x}^T A_{\psi}(\mathbf{\Gamma}^*) \boldsymbol{x}
= 3\sqrt{1-a^2} + \cos(\theta_{i-1,i-2}^*)
\geq 2\sqrt{1-a^2} >0$.
Otherwise $\theta_{i,i+1}^*<\pi/2$, in this case we take $\boldsymbol{x}=[x_1,\ldots,x_n]^T$ such that $x_i=-x_{i+1}=1$ and other entries be $0$,
then $\boldsymbol{x}^T A_{\psi}(\mathbf{\Gamma}^*) \boldsymbol{x}
= 3\sqrt{1-a^2} + \cos(\theta_{i+1,i+2}^*)\geq 2\sqrt{1-a^2}>0$.
Therefore, $A_{\psi}(\mathbf{\Gamma}^*)$ has at least one positive eigenvalue.

Category 2: $\theta_{i,i+1}^*=\alpha$ for any $i\in\mathcal{V}$.
Then due to \rev{the fact} $\Gamma_{i+1}^*=\exp(\alpha\widehat{u})\Gamma_i^*$ for any $i\in\mathcal{V}$,
$\Gamma_1=\exp(n\alpha\widehat{u})\Gamma_1$.
Since $\mathbf{\Gamma}^*\not\in\mathcal{M}_o$, this implies
\[
\alpha\in\{2d\pi/n:d=0,1,\ldots,(n-3)/2\}.
\]
Because $A_{\varphi}(\mathbf{\Gamma}^*)$ is a circulant matrix (see \cite{davis} for details), its eigenvalues can be written explicitly as
\[
\lambda_l = 2\left[ \cos(\alpha) - \cos((l-1)2\pi/n) \right],  \quad l=1,2,\ldots,n.
\]
Then
$\lambda_{(n+1)/2}=\lambda_{(n+3)/2}=\cos(\alpha)-\cos(\pi-\pi/n)>0$.

In conclusion, the system matrix of the linearization system around any $\mathbf{\Gamma}^*\in \mathcal{D}\backslash\mathcal{M}_o$ has at least one positive eigenvalue, indicating that $\mathbf{\Gamma}^*$ is unstable.
\qed
\end{pf}

\begin{rem}
When the inter-agent graph $\mathcal{G}$ is an undirected ring,
the Jacobian matrix of the closed-loop system (\ref{eq:closedloop:psi})-(\ref{eq:closedloop:varphi}) about any equilibrium $\mathbf{\Gamma}^*\in\mathcal{M}_o$ has $3$ zero eigenvalues and $2n-3$ negative eigenvalues.
Similar to the antipodal reduced attitude formation case,
any eigenvector of the zero eigenvalue corresponds to a tangent vector of $\mathcal{M}_o$ at $\mathbf{\Gamma}^*$.
The difference lies in that
%a segment joining a pair of antipodal points has two degree of rotational freedom while a great circle has three degree of rotational freedom.
%
when the number of agents is odd, the structure of $\Gamma_1^*,\ldots,\Gamma_n^*$ has three degrees of rotational freedom in $\mathcal{S}^2$.
\rev{Explicitly, rotating $\Gamma_1^*,\ldots,\Gamma_n^*$ about any axis, which can be either orthogonal or coplanar to the great circle, through any angle will reach another point in the manifold $\mathcal{M}_o$, that is, the tangent space of any point on $\mathcal{M}_o$ has three degrees of freedom. }
\end{rem}

%%%%%%%%%%%%%%%%%%%%%%%%%%%%%%%%%%%%%%%%%%%%%%%%%%%%%%%%%%%%%%%%%%%%%%%%%%%%%%%%
\section{Simulations}\label{sec:simulation}

In this section, we present numerical examples that demonstrate the convergence of the proposed control law (\ref{eq:controller}).

\begin{figure}%[thpb]
  \mbox{}\hfill
  \subfigure[Undirected ring graph]{%, bb=0 0 299 300
    \includegraphics[scale=0.35]{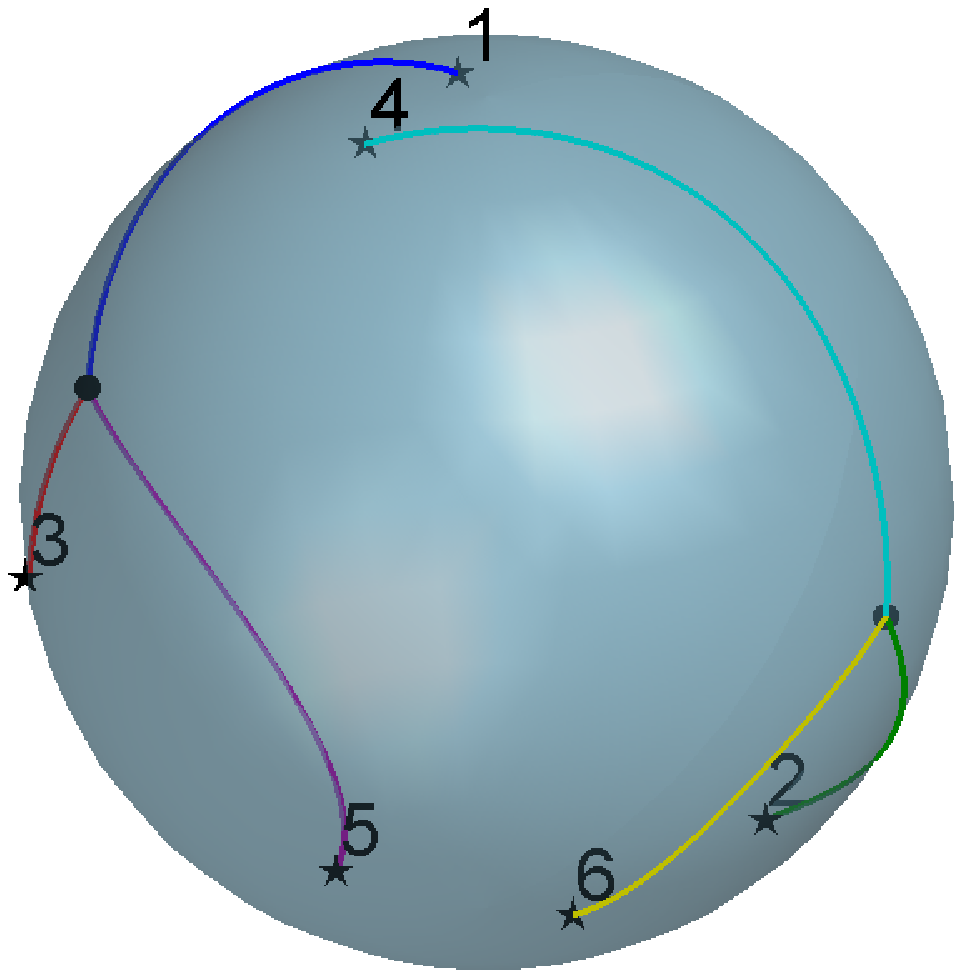}%
    \label{fig:even:a}}
  \hfill
  \subfigure[Directed ring graph]{%
    \includegraphics[scale=0.35]{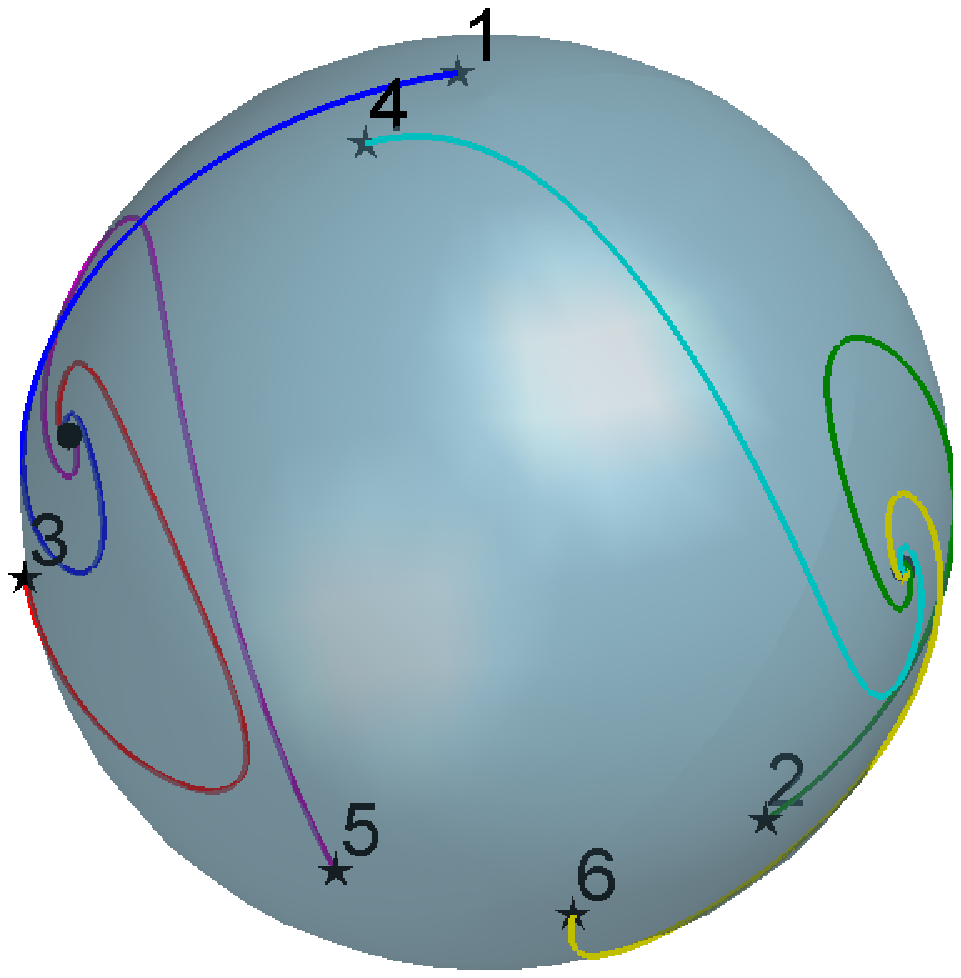}%
    \label{fig:even:b}}
  \hfill\mbox{}
  \caption{Two antipodal reduced attitude formation trajectories for a system of $n=6$ agents. The initial reduced attitudes are marked by pentagrams and the final states are marked by circles.}
  \label{fig:even}
\end{figure}

Concerning the antipodal reduced attitude formation, we consider a system of $6$ rigid-body agents.
The resulting trajectories are illustrated in Fig.~\ref{fig:even},
where Fig.~\ref{fig:even:a} and Fig.~\ref{fig:even:b} share identical
initial states but with the undirected and directed ring inter-agent graph, respectively.
In both cases, the reduced attitudes converge to a pair of antipodal points on $\mathcal{S}^2$, one with oddly indexed agents and the other with evenly indexed agents.

%\begin{figure}
%\begin{center}
%\includegraphics[height=4cm,bb=0 0 232 373]{jcaesar.eps}    % The printed column
%\caption{Gaius Julius Caesar, 100--44 B.C.}  % width is 8.4 cm.
%\label{fig1}                                 % Size the figures
%\end{center}                                 % accordingly.
%\end{figure}
%
%% OR
%
%%\begin{figure}
%%\begin{center}
%%\epsfig{file=jcaesar,width=7cm}
%%\caption{Gaius Julius Caesar, 100--44 B.C.}
%%\label{fig1}
%%\end{center}
%%\end{figure}

\begin{figure}%[thpb]
  \mbox{}\hfill
  \subfigure[$t=0$ (sec)]{%, bb=0 0 299 300
    \includegraphics[scale=0.35]{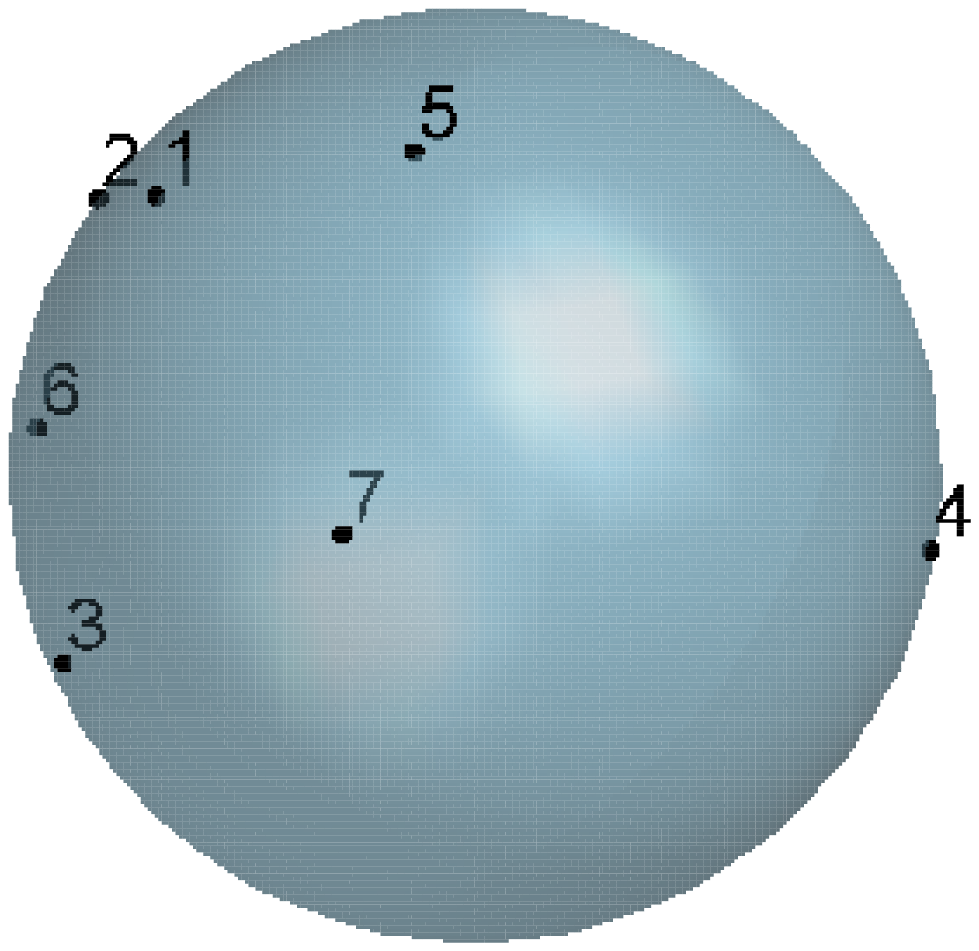}%
    \label{fig:odd:a}}
  \hfill
  \subfigure[$t=40$ (sec)]{%
    \includegraphics[scale=0.35]{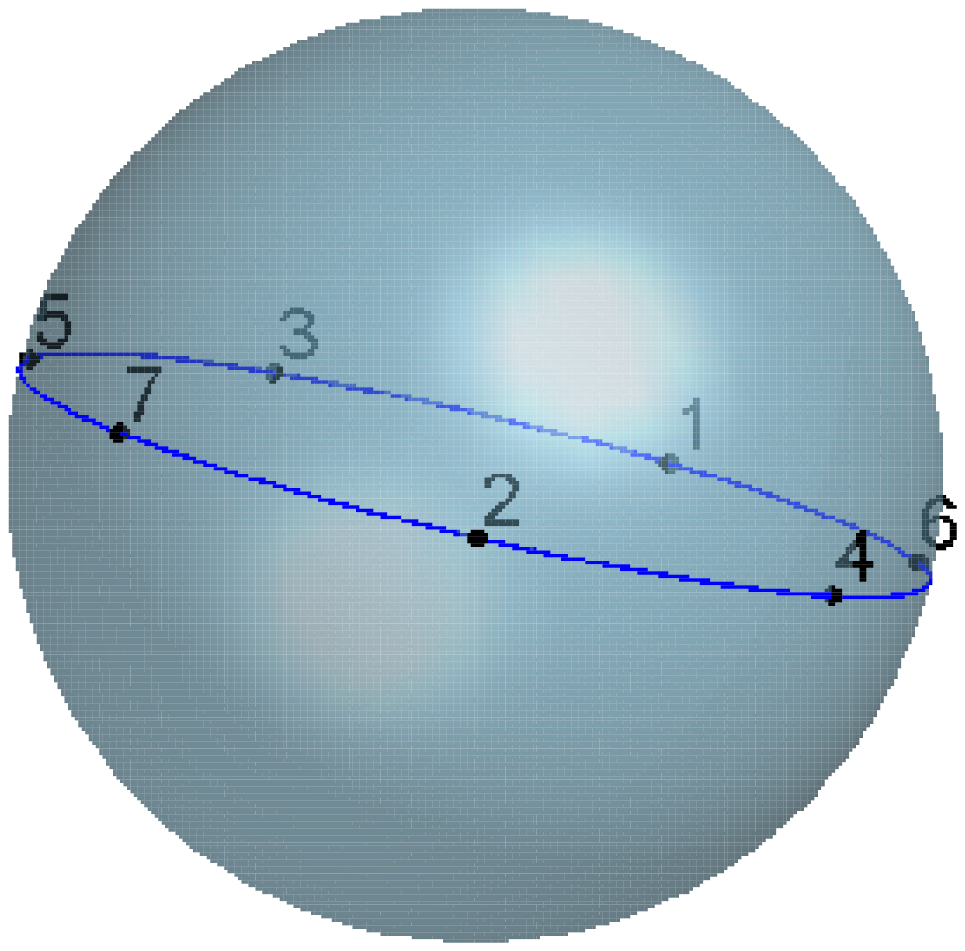}%
    \label{fig:odd:b}}
  \hfill\mbox{}\newline
  \mbox{}\hfill
  \subfigure[$t=45$ (sec)]{%
    \includegraphics[scale=0.35]{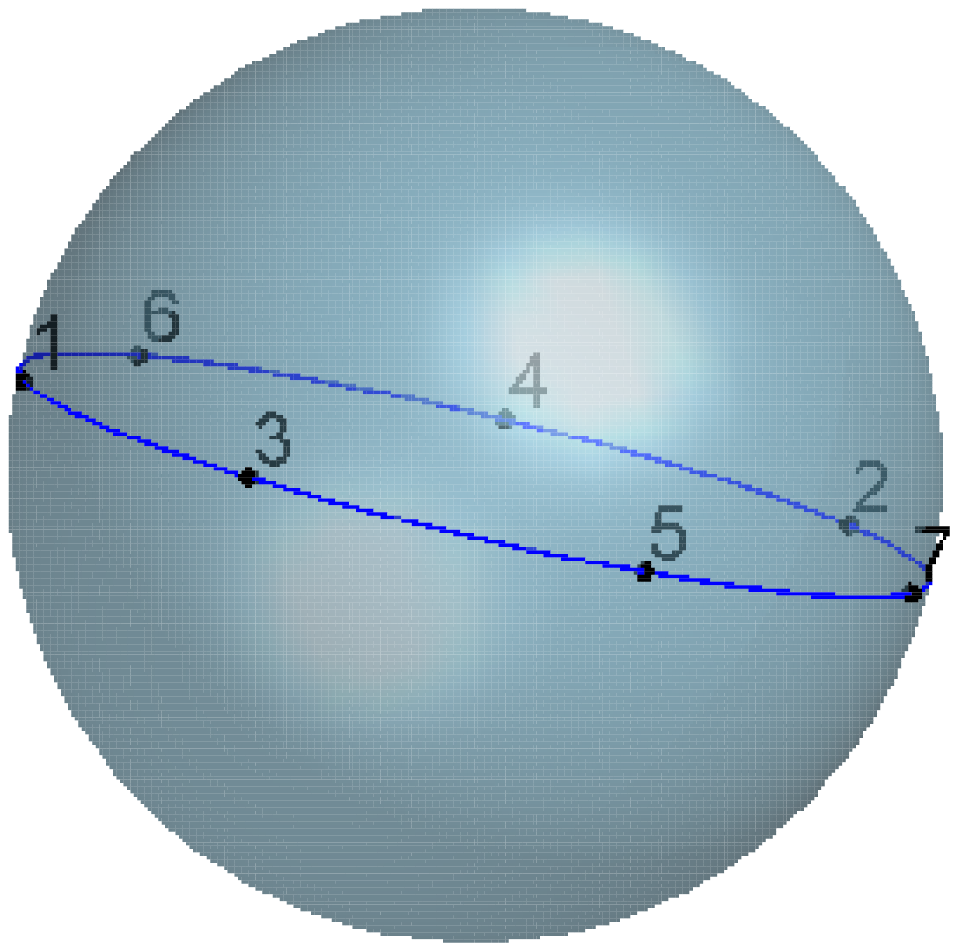}%
    \label{fig:odd:c}}
  \hfill
  \subfigure[$t=50$ (sec)]{%
    \includegraphics[scale=0.35]{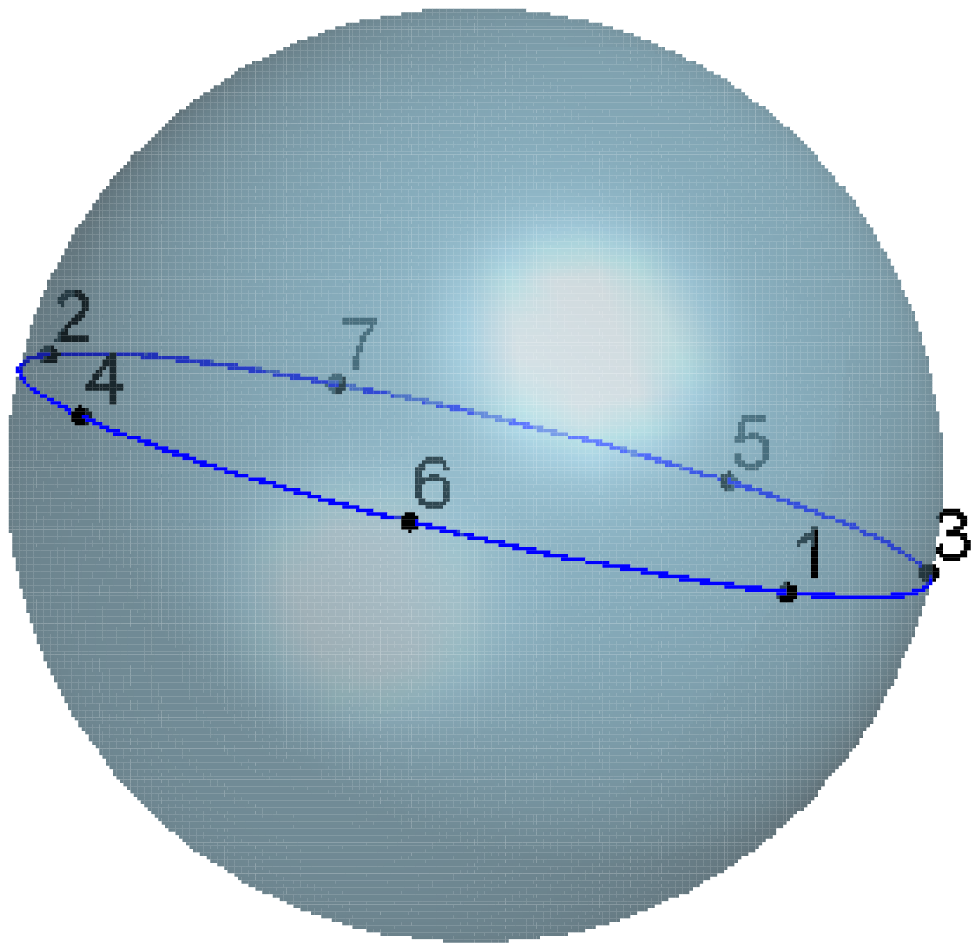}%
    \label{fig:odd:d}}
  \hfill\mbox{}\newline

    \subfigure[Norms of angular velocities]{%
    \includegraphics[scale=0.6]{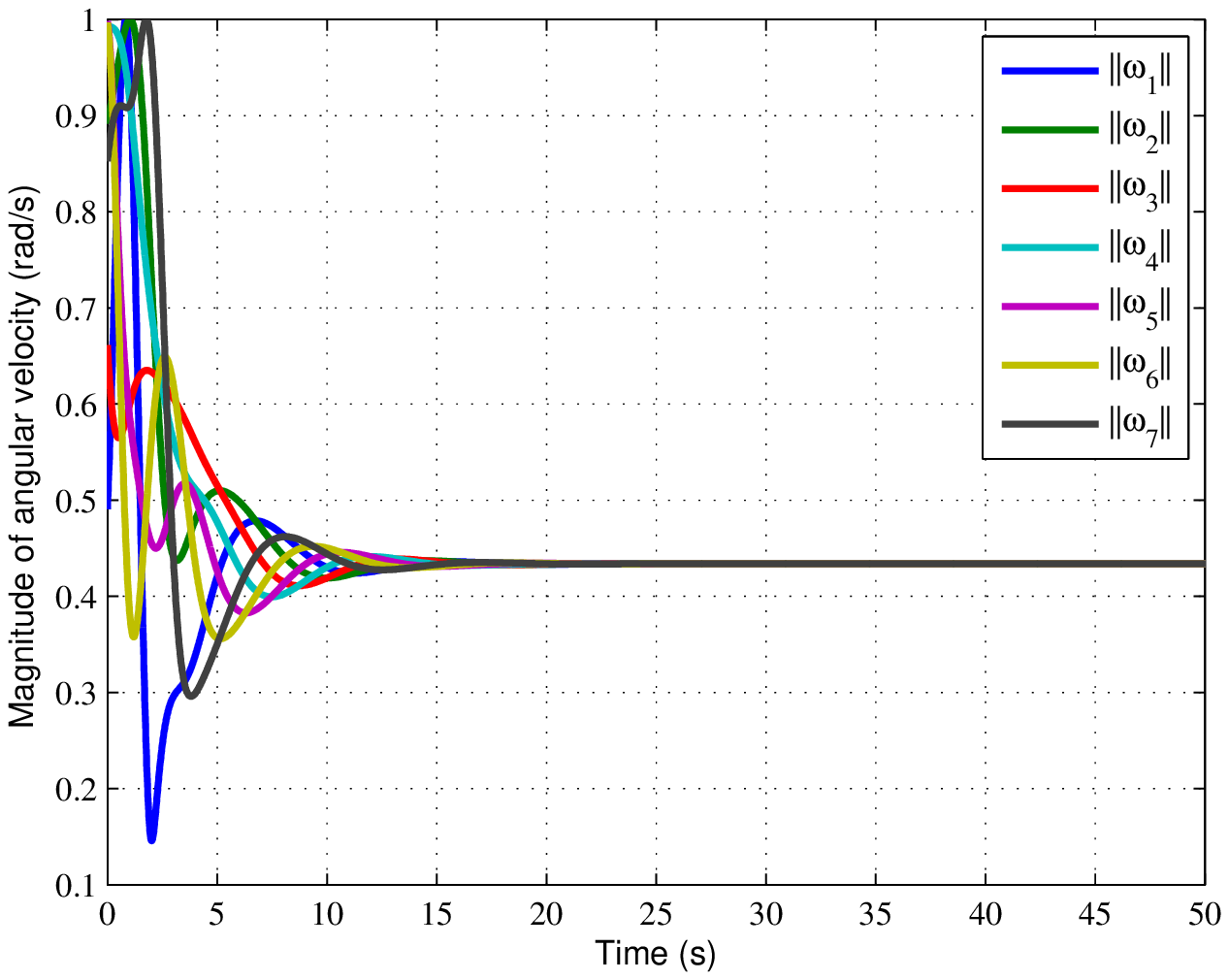}%
    \label{fig:odd:e}}

  \hfill\mbox{}
  \caption{Cyclic reduced attitude formation for a system of $n=7$ agents with the directed ring inter-agent graph. (a)-(d) are the trajectories of reduced attitudes, and (e) is the trajectories of the magnitude of their angular velocities.}
  \label{fig:odd}
\end{figure}

Concerning the cyclic reduced attitude formation, we consider a system
of $7$ rigid-body agents with the directed ring inter-agent graph.
The reduced attitudes initially randomly distribute on $\mathcal{S}^2$ as illustrated in Fig.~\ref{fig:odd:a} and converge to a great circle with equidistant distribution as shown in Fig.~\ref{fig:odd:b}-\ref{fig:odd:d}.
%The cyclic formation is rotating about the axis of the great circle with the period as $2\pi/\sin(\frac{\pi}{7})\approx 14.5$ (sec).
\rev{Fig.~\ref{fig:odd:e} illustrates the time response of the magnitudes of the angular velocities of the agents, which validates that the cyclic formation is rotating about the axis of the great circle with the angular velocity as $\sin(\frac{\pi}{7})\approx 0.44$ (rad/s).}

%%%%%%%%%%%%%%%%%%%%%%%%%%%%%%%%%%%%%%%%%%%%%%%%%%%%%%%%%%%%%%%%%%%%%%%%%%%%%%%%
\section{Conclusions}\label{sec:conclusion}

This paper investigates the reduced attitude formation problem for a system of rigid-body agents with the ring inter-agent graph.
A distributed control law is proposed to stabilize the formation in which
the reduced attitudes of any neighboring agents are at the furthest geodesic distance form one another.
When the total number of agents is even, the agents are
divided into two groups with antipodal reduced attitudes.
When the total number of agents is odd, the reduced attitudes of the
system distribute equidistantly on a great circle, and a rotating
cyclic formation is achieved provided that the ring inter-agent graph
is directed.

As the inter-agent topology varies from a ring, other interesting formations may arise.
Our recent preliminary work shows that for a more general formation, the communication graph should be selected from the class of
graphs sharing the same symmetries with the formation, which suggests
that the more symmetries a formation possesses, the easier it can be
achieved by our approach. Moreover, as the existence of "intrinsic"
formations is partially caused by the compactness of $\mathcal{S}^2$,
more formations could be obtained in other compact spaces following
this "intrinsic" idea.

%% Place acknowledgements here.
%\begin{ack}
%The authors gratefully acknowledge the contribution of National Research Organization and reviewers' comments.
%\end{ack}

%%%%%%%%%%%%%%%%%%%%%%%%%%%%%%%%%%%%%%%%%%%%%%%%%%%%%%%%%%%%%%%%%%%%%%%%%%%%%%%%

\bibliographystyle{plain}        % Include this if you use bibtex %agsm
\bibliography{ref}           % and a bib file to produce the
                                 % bibliography (preferred). The
                                 % correct style is generated by
                                 % Elsevier at the time of printing.

%\begin{thebibliography}{99}     % Otherwise use the
                                 % thebibliography environment.
                                 % Insert the full references here.
                                 % See a recent issue of Automatica
                                 % for the style.
%  \bibitem[Heritage, 1992]{Heritage:92}
%     (1992) {\it The American Heritage.
%     Dictionary of the American Language.}
%     Houghton Mifflin Company.
%  \bibitem[Able, 1956]{Abl:56}
%     B.~C.~Able (1956). Nucleic acid content of macroscope.
%     {\it Nature 2}, 7--9.
%  \bibitem[Able {\em et al.}, 1954]{AbTaRu:54}
%     B.~C. Able, R.~A. Tagg, and M.~Rush (1954).
%     Enzyme-catalyzed cellular transanimations.
%     In A.~F.~Round, editor,
%     {\it Advances in Enzymology Vol. 2} (125--247).
%     New York, Academic Press.
%  \bibitem[R.~Keohane, 1958]{Keo:58}
%     R.~Keohane (1958).
%     {\it Power and Interdependence:
%     World Politics in Transition.}
%     Boston, Little, Brown \& Co.
%  \bibitem[Powers, 1985]{Pow:85}
%     T.~Powers (1985).
%     Is there a way out?
%     {\it Harpers, June 1985}, 35--47.

%\end{thebibliography}

\appendix
% Each appendix must have a short title.
% Sections and subsections are supported in the appendices.

\section{Proofs}\label{appen:sec:proof}

\subsection{Proof of Lemma~\ref{lem:basic:even}}
Due to \rev{the fact}  $W(\mathbf{\Gamma})\leq \pi$ for any $\mathbf{\Gamma} \in (\mathcal{S}^2)^n$, we only need to show $W(\mathbf{\Gamma})=\pi$ if and only if $\mathbf{\Gamma}\in\mathcal{M}_e$.

(Sufficiency)
Take any $\mathbf{\Gamma}\in\mathcal{M}_e$, then $\Gamma_i=(-1)^{i-1}\Gamma_1$ for any $i\in\mathcal{V}$.
Since $n$ is even, $\Gamma_n=-\Gamma_1$.
Hence, $d_{\mathcal{S}^2} (\Gamma_i,\Gamma_{i+1}) = d_{\mathcal{S}^2}(\Gamma_1,-\Gamma_1)=\pi$ for any $i\in\mathcal{V}$, which implies $W(\mathbf{\Gamma})=\pi$.

(Necessity)
Suppose $W(\mathbf{\Gamma})=\pi$, then $d_{\mathcal{S}^2}(\Gamma_i,\Gamma_{i+1})=\pi$ for any $i\in\mathcal{V}$. Hence, $\Gamma_i=-\Gamma_{i+1}$ for any $i\in\mathcal{V}$, which indicates $\mathbf{\Gamma}\in\mathcal{M}_e$.

\subsection{Proof of Lemma~\ref{lem:basic:odd}}

Take any $\mathbf{\Gamma}\in(\mathcal{S}^2)^n$. We first show that $W(\mathbf{\Gamma})\leq \pi-\pi/n$.
\rev{
By the triangular inequality $\theta_{n1} \leq \theta_{13} + \theta_{35} + \cdots + \theta_{n-2,n}$
and
$\theta_{2i-1,2i} + \theta_{2i,2i+1} + \theta_{2i-1,2i+1} \leq 2\pi, i=1,2,\ldots,(n-1)/2$ which followed from Lemma~\ref{lem:sphere}, we obtain
}
%This implies
%\[
%\sum_{i=1}^{n-1} \theta_{i,i+1} + \sum_{i=1}^{(n-1)/2} \theta_{2i-1,2i+1} \leq (n-1)\pi.
%\end{equation}
%On the other hand, by triangle inequality of metric, we get
%\[
%\theta_{1n} \leq \sum_{i=1}^{(n-1)/2} \theta_{2i-1,2i+1}.
%\end{equation}
\begin{eqnarray}
\label{eq:temp0}
\sum_{i\in\mathcal{V}} \theta_{i,i+1}
&\leq& \sum_{i=1}^{n-1} \theta_{i,i+1} + \sum_{i=1}^{(n-1)/2} \theta_{2i-1,2i+1}
\\
\nonumber
&=& \sum_{i=1}^{(n-1)/2} \theta_{2i-1,2i} + \theta_{2i,2i+1} + \theta_{2i-1,2i+1}
\\
\label{eq:temp1}
&\leq& (n-1)\pi,
\end{eqnarray}
which implies $W(\mathbf{\Gamma}) = \min_{i\in\mathcal{V}}\theta_{i,i+1} \leq \pi-\pi/n$.

Next we show $W(\mathbf{\Gamma})=\pi-\pi/n$ if and only if $\mathbf{\Gamma}\in\mathcal{M}_o$.

(Sufficiency)
Take any $\mathbf{\Gamma}\in\mathcal{M}_o$, then
\[
\Gamma_i = \exp\left((i-1)(\pi-\pi/n)\widehat{u}\right) \Gamma_1, \quad\forall i\in\mathcal{V}
\]
for some $u\in\mathcal{S}^2$ and $u^T\Gamma_1=0$.
Because $n$ is odd,
$\Gamma_1 = \exp\left((n-1)\pi\widehat{u}\right)\Gamma_1
= \exp\left((\pi-\pi/n)\widehat{u}\right)\Gamma_n$.
Then,
\begin{eqnarray}
d_{\mathcal{S}^2} (\Gamma_i,\Gamma_{i+1})
&=& d_{\mathcal{S}^2}\left(\Gamma_i,\exp\left((\pi-\pi/n)\widehat{u}\right)\Gamma_i\right)
\nonumber \\
&=& \pi-\pi/n, \quad \forall i\in\mathcal{V}.\nonumber
\end{eqnarray}
Hence, $W(\mathbf{\Gamma})=\pi-\pi/n$.

(Necessity)
Suppose $W(\mathbf{\Gamma})=\pi-\pi/n$.
Since
\[
\sum_{i\in\mathcal{V}} \bar{\theta}_{i,i+1} \leq (n-1)\pi,  \quad\forall \bar{\mathbf{\Gamma}}\in(\mathcal{S}^2)^n,
\]
it follows that $\theta_{i,i+1}=\pi-\pi/n$ for any $i\in\mathcal{V}$.
Then due to \rev{the fact} $\sum_{i\in\mathcal{V}} \theta_{i,i+1} = (n-1)\pi$, both equal signs in (\ref{eq:temp0}) and (\ref{eq:temp1}) hold.
This implies $\theta_{2i-1,2i+1}=2\pi/n$ for any $i\in\{1,2,\ldots,(n-1)/2\}$ and $\Gamma_1,\Gamma_2,\ldots,\Gamma_n$ lie on a great circle of $\mathcal{S}^2$.
%by Lemma~\ref{lem:sphere}.
Let $u=-\widehat{\Gamma}_1\Gamma_n/\sin(\pi/n)$.
Because $\theta_{1n}=\sum_{i=1}^{(n-1)/2}\theta_{2i-1,2i+1}=\pi-\pi/n$, $\Gamma_1,\Gamma_3,\Gamma_5,\ldots,\Gamma_n$ lie equidistantly on the shorter arc of the great circle joining $\Gamma_1$ and $\Gamma_n$. Therefore,
\[
\Gamma_{2i+1}=\exp\left(2(\pi-\pi/2)\widehat{u}\right) \Gamma_{2i-1}
\]
holds for $i=1,2,\ldots,(n-1)/2$.
Because $\theta_{2i-1,2i} + \theta_{2i,2i+1} + \theta_{2i-1,2i+1} = 2\pi$, $\Gamma_{2i}$ is the middle point of the longer arc of the great circle joining $\Gamma_{2i-1}$ and $\Gamma_{2i+1}$.
Therefore,
\[
\Gamma_{2i}=\exp\left((\pi-\pi/n)\widehat{u}\right) \Gamma_{2i-1}
\]
holds for $i=1,2,\ldots,(n-1)/2$.
As a result, $\Gamma_i=\exp\left((i-1)(\pi-\pi/n)\widehat{u}\right) \Gamma_1$
for any $i\in\mathcal{V}$, i.e., $\mathbf{\Gamma}\in\mathcal{M}_o$.

\subsection{Proof of Lemma~\ref{lem:odd:bound:0}}
Let
$\bar{\Gamma}_i= \exp((i-1)(\pi-\pi/n)\widehat{k}_{12}) \Gamma_{1}$
for any $i\in\mathcal{V}$.
Then $\bar{\mathbf{\Gamma}}=\{\bar{\Gamma}_i\}_{i\in\mathcal{V}}\in \mathcal{M}_o$ and we will show that $d_{(\mathcal{S}^2)^n}(\mathbf{\Gamma},\bar{\mathbf{\Gamma}})\leq 2n\nu$.

Clearly,
$d_{\mathcal{S}^2}(\Gamma_{1},\bar{\Gamma}_{1})=0$ and
\begin{eqnarray}
d_{\mathcal{S}^2}(\Gamma_{2},\bar{\Gamma}_{2})
&=& d_{\mathcal{S}^2}(\exp(\theta_{12}\widehat{k}_{12})\Gamma_1,\exp((\pi-\pi/n)\widehat{k}_{12}) \Gamma_{1})
\nonumber \\
&=&|\pi-\pi/n-\theta_{12}| \leq \nu^2. \nonumber
\end{eqnarray}
\rev{
To estimate the upper bound of $d_{\mathcal{S}^2} (\Gamma_i,\bar{\Gamma}_{i})$ for $i\geq 3$, we divide the remaining agents into two groups $3\leq i\leq (n+1)/2$ and $(n+3)/2 \leq i \leq n$, and then consider each group respectively.
In the following,}
we denote $c_{\beta}$ and $s_{\beta}$ as the respective abbreviations for $\cos(\beta)$ and $\sin(\beta)$.

For agent $i$ in the first group, that is, $3\leq i\leq (n+1)/2$, we have that
\revA{
\begin{eqnarray}
&&{}d_{\mathcal{S}^2} (\Gamma_i,\bar{\Gamma}_{i})
\nonumber \\
&=& d_{\mathcal{S}^2}(\exp(\theta_{1i}\widehat{k}_{1i})\Gamma_1,
       \exp((i-1)(\pi-\pi/n)\widehat{k}_{12})\Gamma_{1})
\nonumber \\
&\leq& d_{\mathcal{S}^2}(\exp(\theta_{1i}\widehat{k}_{1i})\Gamma_1,
       \exp((i-1)\theta_{12}\widehat{k}_{12})\Gamma_{1})
\nonumber\\
&& \, +  d_{\mathcal{S}^2}(\exp((i-1)\theta_{12}\widehat{k}_{12})\Gamma_{1},
        \exp((i-1)(\pi-\pi/n)\widehat{k}_{12})\Gamma_{1})
\nonumber\\
&\leq&  d_{\mathcal{S}^2}(\exp(\theta_{1i}\widehat{k}_{1i})\Gamma_1,
       \exp((i-1)\theta_{12}\widehat{k}_{12})\Gamma_{1})
\nonumber\\
&& \quad +  (i-1)|\pi-\pi/n-\theta_{12}|.
\label{eq:2016022801}
\end{eqnarray}
}
By applying Lemma~\ref{lem:sphere} \revA{to} the points
$\exp(\theta_{1i}\widehat{k}_{1i})\Gamma_1$,
$\exp((i-1)\theta_{12}\widehat{k}_{12})\Gamma_1$,
$\Gamma_1\in\mathcal{S}^2$, we obtain
\[
d_{S^2}(\exp (\theta_{1i}\widehat{k}_{1i})\Gamma_1, \exp((i-1)\theta_{12}\widehat{k}_{12})\Gamma_1)=\arccos(f_i),
\]
where $f_i=c_{\theta_{1i}}c_{(i-1)\theta_{12}} + s_{\theta_{1i}}s_{(i-1)\theta_{12}}k_{1i}^T k_{12}$. Then (\ref{eq:2016022801}) can be rewritten as
\begin{equation}\label{eq:2016022802}
d_{\mathcal{S}^2} (\Gamma_i,\bar{\Gamma}_{i})
\leq  \arccos(f_i) +  (i-1)|\pi-\pi/n-\theta_{12}|.
\end{equation}
We next estimate a lower bound of $f_i$.
By Lemma~\ref{lem:sphere},
\begin{equation}\label{eq:temp3}
c_{\theta_{2i}}
= c_{\theta_{12}}c_{\theta_{1i}} + s_{\theta_{12}}s_{\theta_{1i}}k_{1i}^Tk_{12},
\end{equation}
therefore,
$
f_i= - c_{\theta_{1i}} s_{(i-2)\theta_{12}}/s_{\theta_{12}}
+c_{\theta_{2i}} s_{(i-1)\theta_{12}}/s_{\theta_{12}}.
$

\revA{Next, we discuss the lower bound of $f_i$ for the case when $i$ is odd and even, respectively.}

If $i$ is odd, by Lemma~\ref{lem:sphere} and (\ref{eq:odd:temp1}), %it can be showed that
\[
 - (n-i+1)\nu^2 \leq \theta_{1i}-(i-1)\pi/n \leq  (i-1)\nu^2.
\]
\revA{
Combine this and the inequality
\[
 \pi-\frac{1}{n}\pi-\nu^2 \leq \theta_{i,i+1} \leq  \pi-\frac{1}{n}\pi+\nu^2,
\]
which comes from (\ref{eq:odd:temp1}).
For the group $3\leq i\leq (n+1)/2$, the lower bound of $f_i$ can be analyzed as follows.
}

(i) When $i=3$,
\revA{
\begin{eqnarray}
f_i
&&= - c_{\theta_{13}} + 2 c_{\theta_{12}} c_{\theta_{23}}
\geq - c_{\frac{2}{n}\pi-(n-2)\nu^2} + 2 c_{\pi-\frac{1}{n}\pi-\nu^2}^2
%\\
%&= - c_{\frac{2}{n}\pi-(n-2)\nu} + 2 c_{\frac{1}{n}\pi+\nu}^2
%\\ &= 1 - c_{\frac{2}{n}\pi-(n-2)\nu^2} + c_{\frac{2}{n}\pi+2\nu^2}
\nonumber \\
&&=- c_{\frac{2}{n}\pi-(n-2)\nu^2} + 1 + c_{\frac{2}{n}\pi+2\nu^2}
= 1 - 2s_{\frac{n}{2}\nu^2} s_{\frac{2}{n}\pi-\frac{n-4}{2}\nu^2}
\nonumber \\
&&\geq 1- 2s_{\frac{n}{2}\nu^2}
\geq 1- n\nu^2 \quad \quad \quad \quad(since \, x \geq s_x, \forall x \geq 0)
\nonumber \\
&&\geq c_{2\sqrt{n}\nu}. \quad \quad \quad \quad \quad \quad\quad \quad \quad \quad (since \, \nu \in [0,\sqrt{2}/n]) \nonumber
\end{eqnarray}
}
(ii) When $5\leq i\leq (n+1)/2$,
\revA{
by the trigonometric identity
\begin{equation}\label{eq:2016093001}
s_{k\alpha}=
\left\{
  \begin{array}{ll}
    s_\alpha [1+\sum_{j=1}^{(k-1)/2} 2c_{2j\alpha}], \quad & \hbox{$k$\, is  odd;} \\
    2s_\alpha  \sum_{j=1}^{k/2} 2c_{(2j-1)\alpha}, & \hbox{$k$\, is  even,}
  \end{array}
\right.
\end{equation}
we have
}
\[%\begin{eqnarray}\label{eq:odd:temp2}
f_i
= - c_{\theta_{1i}} - 2 \sum_{j=1}^{(i-3)/2} c_{\theta_{1i}} c_{2j\theta_{12}}
+ 2\sum_{j=1}^{(i-1)/2} c_{\theta_{2i}} c_{(2j-1)\theta_{12}}.
\]%\end{eqnarray}
By Lemma~\ref{lem:sphere} and (\ref{eq:odd:temp1}), %it can be showed that
\[
 - (i-2)\nu^2 \leq
\theta_{2i}-\pi+ (i-2)\pi/n
\leq  (n-i+2)\nu^2.
\]
Then by monotonicity of $f_i$ with respect to $\theta_{1i}$, $\theta_{2i}$ and $\theta_{12}$, %$f_i$ is larger than the value of (\ref{eq:odd:temp2}) when $\theta_{1i}=\frac{i-1}{n}\pi - (n-i+1)\nu^2$, $\theta_{2i}=\pi-\frac{i-2}{n}\pi - (i-2)\nu^2$, $\theta_{12}=\pi-\frac{1}{n}\pi+\nu^2$ in the second term of (\ref{eq:odd:temp2}) and $\theta_{12}=\pi-\frac{1}{n}\pi-\nu^2$ in the third term of (\ref{eq:odd:temp2}). Simplifying the value yields
we obtain
\revA{
\begin{eqnarray}
f_i
%&=& - c_{\theta_{1i}} -  \sum_{j=1}^{(i-3)/2} (c_{2j\theta_{12}+\theta_{1i}}+c_{2j\theta_{12}-\theta_{1i}})
%\nonumber \\
%&&{}+ \sum_{j=1}^{(i-1)/2} (c_{(2j-1)\theta_{12}+\theta_{2i}} + c_{(2j-1)\theta_{12}-\theta_{2i}} ).
%\nonumber \\
&\geq& - c_{\frac{i-1}{n}\pi - (n-i+1)\nu^2} - \sum_{j=1}^{(i-3)/2} c_{\frac{i-1-2j}{n}\pi - (n-i+1-2j)\nu^2}
\nonumber \\
&&{}
\!\!\!-\!\!\sum_{j=1}^{\!\!\!(i-3)/2} c_{\frac{i-1+2j}{n}\pi - (n-i+1+2j)\nu^2}
+\!\!\sum_{j=1}^{\!\!\!(i-1)/2} c_{\frac{i-3+2j}{n}\pi + (i-3+2j)\nu^2}
\nonumber \\
&&{}
+ \sum_{j=1}^{(i-1)/2} c_{\frac{i-1-2j}{n}\pi + (i-1-2j)\nu^2}
\nonumber \\
&=& 1 - 2 \sum_{j=1}^{(i-3)/2}s_{\frac{n-4j}{2}\nu^2}s_{\frac{i-1-2j}{n}\pi-\frac{n+2-2i}{2}\nu^2}
\nonumber \\
&&{} - 2 \sum_{j=1}^{(i-1)/2}s_{\frac{n+4j-4}{2}\nu^2}s_{\frac{i-3+2j}{n}\pi-\frac{n+2-2i}{2}\nu^2}
\nonumber\\
&\geq& 1 - \sum_{j=1}^{(i-3)/2}(n-4j)\nu^2 - \sum_{j=1}^{(i-1)/2}(n+4j-4)\nu^2
\nonumber\\
&=& 1-(i-2)n\nu^2 \geq c_{2\sqrt{n(i-2)}\nu}\nonumber .
\end{eqnarray}
}
If $i$ is even, \revA{by identity (\ref{eq:2016093001}), we obtain}
\begin{eqnarray}
f_i
&&= - 2 \sum_{j=1}^{(i-2)/2} c_{\theta_{1i}}c_{(2j-1)\theta_{12}}
+ c_{\theta_{2i}} + 2 \sum_{j=1}^{(i-2)/2} c_{\theta_{2i}}c_{2j\theta_{12}}
\nonumber \\
&&\geq c_{2\sqrt{n(i-2)}\nu}. \nonumber
\end{eqnarray}
\rev{
Therefore, combining (\ref{eq:2016022802}) yields
\begin{equation}\label{eq:2016092904}
d_{\mathcal{S}^2} (\Gamma_i,\bar{\Gamma}_{i})
\leq 2\sqrt{n(i-2)}\nu +  (i-1)\nu^2, \nonumber
\end{equation}
which holds for any agent $i$ in the first group.}

\rev{
Next we consider agent $i$ in another group, that is $(n+3)/2\leq i\leq n$.
In this case, because}
$\bar{\Gamma}_{i}
= \exp(-n(\pi-\pi/n)\widehat{k}_{12})\bar{\Gamma}_{i}
=\exp(-(n-i+1)(\pi-\pi/n)\widehat{k}_{12})\Gamma_{1}
$, we obtain
\begin{eqnarray}\label{eq:2016092902}
&&d_{\mathcal{S}^2} (\Gamma_i,\bar{\Gamma}_{i})
\nonumber \\
&&= d_{\mathcal{S}^2}(\exp(\theta_{1i}\widehat{k}_{1i})\Gamma_1, \nonumber\\
   && \qquad\qquad   \exp(-(n-i+1)(\pi-\pi/n)\widehat{k}_{12})\Gamma_{1})
%%%%%%%%%%%%%%%%%%%%%%%%
%%%%%%%%%%%%%%%%%%%%%%%%
\nonumber \\
&&\leq  d_{\mathcal{S}^2}(\exp(\theta_{1i}\widehat{k}_{1i})\Gamma_1,
       \exp(-(n-i+1)\theta_{12}\widehat{k}_{12})\Gamma_{1})
\nonumber \\
&&\qquad\qquad + (n-i+1)|\pi-\pi/n-\theta_{12}|
%%%%%%%%%%%%%%%%%%%%%%%%
%%%%%%%%%%%%%%%%%%%%%%%%
\nonumber \\
&&=  \arccos(f_i) + (n-i+1)|\pi-\pi/n-\theta_{12}|,
\end{eqnarray}

where
$f_i=c_{\theta_{1i}}c_{(n-i+1)\theta_{12}} - s_{\theta_{1i}}s_{(n-i+1)\theta_{12}}k_{1i}^T k_{12}$.
Then combining (\ref{eq:temp3}) and applying similar computations as the first group yields
\begin{eqnarray}
f_i
&=& - c_{\theta_{2i}} s_{(n-i+1)\theta_{12}}/s_{\theta_{12}}
+ c_{\theta_{1i}} s_{(n-i+2)\theta_{12}}/s_{\theta_{12}}.
\nonumber \\
&\geq& c_{2\sqrt{n(n-i+1)}\nu}. \nonumber
\end{eqnarray}
\revA{
Notice that in the above derivation, the inequality
\begin{equation}\label{eq:2016092901}
\sqrt{n(n-1)/2}\nu <\nu n/\sqrt{2}\leq 1
\end{equation}
is used. Therefore, (\ref{eq:2016092902}) implies
}
\begin{equation}\label{eq:2016092903}
d_{\mathcal{S}^2} (\Gamma_i,\bar{\Gamma}_{i})
\leq  2\sqrt{n(n-i+1)}\nu +  (n-i+1)\nu^2
\end{equation}
holds for any agent in the second group.

\revA{
In conclusion, by (\ref{eq:2016092904}), (\ref{eq:2016092901}), (\ref{eq:2016092903}) we obtain}
\begin{eqnarray}
d_{\mathcal{M}_o}(\mathbf{\Gamma})
&\leq& d_{(\mathcal{S}^2)^n}(\mathbf{\Gamma},\bar{\mathbf{\Gamma}})
= \max_{i\in\mathcal{V}} d_{\mathcal{S}^2} (\Gamma_i,\bar{\Gamma}_{i})
\nonumber \\
&\leq& 2\sqrt{n(n-1)/2}\nu +  (n-1)\nu^2/2
\nonumber \\
&\leq& \left(\sqrt{2}n+n/2\right)\nu \leq 2n\nu. \nonumber
\end{eqnarray}
The proof is completed.

%\addtolength{\textheight}{-11cm}
% This command serves to balance the column lengths on the last page of the document manually.
% It shortens the textheight of the last page by a suitable amount.
% This command does not take effect until the next page so it should come on the page before the last.
% Make sure that you do not shorten the textheight too much.

\end{document}